\input ./style/arxiv-general.cfg
\documentclass[aap,MSNbibl,seceqn,dvips]{arximspdf}
\makeatletter
   \@ifpackageloaded{graphicx}{}{\usepackage{graphicx}}
\makeatother

\doi{10.1214/14-AAP1055} 
\volume{25}
\issue{5}
\pubyear{2015}
\firstpage{2535}
\lastpage{2562}
\docsubty{FLA}

\makeatletter
\def\cal{\mathcal}
\newcommand{\eqref}[1]{(\ref{#1})}

\newtheorem{Theorem}{Theorem}[section]
\newproclaim{Definition}[Theorem]{Definition}
\newtheorem{Proposition}[Theorem]{Proposition}
\newtheorem{proposition}{Proposition}[section]

\newproclaim{Assumption}[Theorem]{Assumption}

\newtheorem{Lemma}[Theorem]{Lemma}

\newproclaim{Remark}[Theorem]{Remark}
\newproclaim{remark}[proposition]{Remark}
\newproclaim{Example}[Theorem]{Example}

\def\E{\mathbb{E}}
\def\F{\mathbb{F}}
\def\G{\mathbb{G}}
\def\M{\mathbb{M}}
\def\P{\mathbb{P}}
\def\R{\mathbb{R}}

\def\Cc{{\cal C}}
\def\Dc{{\cal D}}
\def\Fc{{\cal F}}
\def\Nc{{\cal N}}
\def\Pc{{\cal P}}
\def\Uc{{\cal U}}
\def\Yc{{\cal Y}}
\def\Zc{{\cal Z}}

\def\Pb{{\overline\P}}

\def\Mh{{\widehat M}}

\def\a{\alpha}

\def\Om{\Omega}
\def\Omb{\overline{\Omega}}
\def\eps{\varepsilon}
\def\xb{\mathbf{x}}

\def\Yt{\tilde{Y}}
\def\0{\mathbf{0}}

\def\ah{\widehat{a}}

\def\S{\mathbb{S}}
\def\Fcb{\overline{{\cal F}}}
\def\Fbb{\overline{\mathbb{F}}}

\def\Yct{\widetilde{\Yc}}
\def\Zct{\widetilde{\Zc}}

\def\ut{\tilde{u}}

\def\Yt{{\widetilde Y}}

\def\x{\times}

\def\eps{\varepsilon}

\def\reff#1{(\ref{#1})}

\def\1{\mathbf{1}}

\def\x{\times}
\def\1{\mathbf{1}}

\makeatother

\begin{document}
\begin{frontmatter}

\title{Weak approximation of second-order BSDEs}
\runtitle{Weak approximation of 2BSDEs}

\begin{aug}
\author{\fnms{Dylan} \snm{Possama\"{i}}}\ead[label=e1]{possamai@ceremade.dauphine.fr}
\and
\author{\fnms{Xiaolu} \snm{Tan}\corref{}\ead[label=e2]{tan@ceremade.dauphine.fr}}
\runauthor{D. Possama\"{i} and X. Tan}
\affiliation{University of Paris-Dauphine}
\address{Ceremade\\
Universite Paris-Dauphine\\
Place du Marechal de Lattre de Tassigny\\
75775, Paris Cedex 16\\
France\\
\printead{e1}\\
\phantom{E-mail:\ }\printead*{e2}}
\end{aug}

\received{\smonth{9} \syear{2013}}
\revised{\smonth{5} \syear{2014}}


\begin{abstract}
We study the weak approximation of the second-order backward SDEs (2BSDEs),
when the continuous driving martingales are approximated by discrete
time martingales.
We establish a convergence result for a class of 2BSDEs,
using both robustness properties of BSDEs, as proved in Briand, Delyon
and M\'emin [\textit{Stochastic Process. Appl.}
\textbf{97}
(2002)
229--253], and tightness of solutions to
discrete time BSDEs.
In particular, when the approximating martingales are given by some
particular controlled Markov chains, we obtain several concrete
numerical schemes for 2BSDEs, which we illustrate on specific examples.
\end{abstract}

%

%
\end{frontmatter}

\section{Introduction}\label{sec1}

Weak approximation is an important technique in\break stochastic analysis. A
famous and classical result in this spirit is Donsker's theorem which
stipulates the following. Let $(\zeta_k)_{k \ge1}$ be a sequence of
i.i.d. centered random variables such that $\operatorname{Var}(\zeta
_1) = 1$, and define
\[
S^n_t:= \frac{1}{\sqrt{n}} \sum
_{k=1}^{[nt]} \zeta_k,
\]
then the process $S^n_{\cdot}$ converges weakly to a Brownian motion
$W$. In particular, suppose that $f \dvtx\R\longrightarrow\R$ is a
bounded continuous function, we then have the following convergence:
\[
\E\bigl[ f\bigl(S^n_T\bigr) \bigr] \to\E\bigl[
f(W_T)\bigr].
\]
Similar result have been obtained for diffusion processes defined as
solutions to stochastic differential equations (SDEs in the sequel);
see, for example, Jacod and Shiryaev \cite{JacodShiryaev}.
We also remind the reader that in this Markovian setting, the value
$\E[f(W_T)]$ can be characterized using the heat equation from the
Feynmann--Kac formula.

Backward stochastic differential equations (BSDEs in the sequel),
which were introduced by Pardoux and Peng \cite{PardouxPeng}, as well
as the more recent notion of $G$-expectation of Peng \cite
{Peng_G_expec}, are particular cases of so-called nonlinear
expectations, and their weak approximation properties have attracted a
lot of attention in the recent years.
Hence, in Briand, Delyon and M\'emin \cite{BriandDelyonMemin}, the
authors studied the convergence of the solutions of the BSDE when the
driving Brownian motion is approximated by a sequence of martingales.
In particular, when the Brownian motion is approximated by some random
walks, they obtained a weak convergence result similar to the above
Donsker's theorem.
More recently, Dolinsky, Nutz and Soner \cite{Dolinsky_Nutz_Soner}
studied the weak approximation of $G$-expectation. Since
$G$-expectation can be considered as a sublinear expectation on the
canonical space of continuous trajectories, by the analogue of
Donsker's theorem, they approximated it by a sequence of sublinear
expectations on the canonical space of discrete time paths.
Extending BSDE and $G$-expectation, the second-order backward SDEs
(2BSDEs) introduced by Soner, Touzi and Zhang \cite
{Soner_Touzi_Zhang1}, can be represented as the supremum of a family of
nonlinear expectations on the canonical space of continuous trajectories.
In particular, it generalizes the Feynmann--Kac formula to the fully
nonlinear case.
We are then motivated to extend the weak approximation property to 2BSDEs.

We notice that the weak approximation property should be an important
property of the continuous time dynamic models,
when it is the continuous limit of discrete time models.
For example, in finance, it is convenient to use a Brownian motion to
model the evolution of a risky asset, despite the fact that such a
price only exists on discrete time instants.
Therefore, it is important to confirm that as we take the limit of the
discrete time model, it converges to the continuous time model.

Finally, weak approximation is also an important technique in
numerical analysis;
see, for example, Kushner and Dupuis \cite{KushnerDupuis} in the
context of stochastic control problems,
and Dolinsky \cite{Dolinsky} for pricing the financial ``game'' options.
The main idea is to interpret the numerical scheme as a controlled
Markov chain system,
which converges weakly to the continuous time system.
We notice also that another point of view is from the PDEs,
which characterizes the solution of these dynamic problems in the
Markovian case.
A powerful numerical analysis method in this context is the monotone
convergence theorem of Barles and Souganidis \cite{Barles_Souganidis}.
Comparing to the PDE numerical methods, the weak approximation method
permits usually to relax regularity and integrability conditions, and
also permits to study the non-Markovian problems as shown in Tan \cite
{TanNumStoContr}.

The main contribution of the paper is to prove a weak approximation
property for a class of 2BSDEs, which can be considered as an extension
of Donsker's theorem in this nonlinear context.
Further, using some controlled Markov chains as approximating
martingales, we obtain some numerical schemes for a class of 2BSDEs.
In particular, these numerical schemes are coherent with the classical
schemes proposed for the nonlinear PDEs in the Markovian cases.
We also notice that these related numerical schemes have been largely
tested in the previous literature; see, for example, Fahim, Touzi and
Warin \cite{Fahim_Touzi_Warin}, Tan \cite{TanSplitting}, Guo, Zhang and
Zhuo \cite{GuoZhangZhuo}, etc.

The rest of the paper is organized as follows.
In Section~\ref{sec:main_result}, we introduce the class of 2BSDEs
that is studied in the paper,
and give first an equivalence result using two different classes of
driving martingales. By considering a sequence of discrete time
equations, we give a general weak approximation result, that is, the
discrete time solution converges to the solution of a class of 2BSDE.
Then in Section~\ref{sec:num}, by considering some particular
controlled Markov chains, we can interpret the discrete time equations
as numerical schemes, and the weak approximation result justifies the
convergence of the numerical schemes. Section~\ref{sec:num2} is devoted
to some numerical examples, highlighting the convergence of the
proposed numerical schemes.
In Section~\ref{sec:equiv}, we complete the proof of the equivalence theorem,
and finally in Section~\ref{sec:cvg}, we report the proof of the weak
approximation theorem.

Throughout the paper, we use the following notation. For
every $(x,y)\in\R^d\times\R^d$, we denote by $x\cdot y$ the usual
scalar product of $x$ and $y$, and for any $(x,y)\in\R^{d\x d}\times
\R
^{d\x d}$, we denote by $x \dvtx y:=\operatorname{Tr}(xy)$. Similarly, $x^T$
will denote the usual transposition and $|x|$ the Euclidean norm in
the corresponding space.

\section{The 2BSDE and its weak approximation}
\label{sec:main_result}

In this section, we first introduce the class of second-order BSDEs
that we next propose to approximate by the supremum of a family of
BSDEs driven by approximating discrete time martingales.
A convergence result is given under sufficient conditions, while the
proof is postponed to other sections.

\subsection{A class of 2BSDEs}

Let $\Om:= \{ \omega\in C([0,T],\R^d) \dvtx\omega_0 = 0 \}$
denote the
canonical space of continuous paths on $[0,T]$ which start at $0$,
$B$ be the canonical process, $\F= (\Fc_t)_{0 \le t \le T}$ the
canonical filtration and $\P_0$ the Wiener measure on $\Om$ under which
$B$ is a standard Brownian motion.
Denote by $\F^+ = (\Fc^+_t)_{0 \le t \le T}$ the right-continuous
filtration defined by $\Fc_t^+ := \bigcap_{s > t} \Fc_s$ for all $t < T$
and $\Fc^+_T = \Fc_T$.
For every probability measure $\P$ on $(\Omega,\mathcal F_T)$, we
denote by $\Fbb^{\P}$ the $\P$-augmented filtration of $\F$ and
$\overline{\F^+}^\P$ the $\P$-augmented filtration of $\F^+$. Moreover,
for any $\mathbf{x}\in\Omega$, and for any $t\in[0,T]$, we note
$\|x\|_t:=\sup_{0\leq s\leq t}|x_s|$.
A probability measure $\P$ on $\Omega$ such that $B$ is a $\P$-local
martingale will be called a local martingale measure.

We recall that by results of Bichteler \cite{bichteler} (see also
Karandikar \cite{kar} for a simplified exposition) there are two $\F
$-progressive processes on $\Om$ given by
\[
\langle B \rangle _t := B_t B_t^T
- 2 \int_0^t B_s
\,dB_s^T \quad\mbox{and}\quad \ah_t := \limsup
_{\eps\downarrow0} \frac{1}{\eps} \bigl( \langle B \rangle
_t - \langle B \rangle _{t-\eps} \bigr),
\]
such that $\langle B \rangle $ coincides with the $\P$-quadratic
variation of
$B$, $\P$-a.s., for all local martingale measures $\P$.

We consider next a set $A$ such that
%
\begin{equation}
A \subset\S_d^+ \mbox{ is compact, convex and } a \ge
\eps_0 I_d, \forall a \in A,
\end{equation}
where $\S_d^+$ is the set positive, symmetric $d\x d$ matrices and
where $\eps_0 > 0$ is a fixed constant.
We denote by $\Pc_W$ the collection of all local martingale measures
$\P$ such that $\ah\in A$, $d\P\x \,dt$-a.e.,
and by $\Pc_S \subset\Pc_W$ the subset consisting of all probability measures
\[
\P^{\alpha} := \P_0 \circ\bigl(X^{\alpha}
\bigr)^{-1}\qquad \mbox{where } X^{\alpha}_t := \int
_0^t \alpha^{1/2}_s
\,dB_s, \P_0\mbox{-a.s.}
\]
for some $\F$-progressively measurable process $\alpha$ taking values
in $A$.

Let now $\xi\dvtx\Om\to\R$ be a random variable, $g \dvtx[0,T]
\x\Om
\x\R\x\R^d \x\S_d^+ \to\R^{d\x d}$ be a function which will play
the role of our generator. Then for every $\P\in\Pc_W$, we consider
the following generalized BSDE under $\P$:
%
\begin{eqnarray}
\label{eq:BSDE} \Yc_t^\P &=& \xi(B_{\cdot}) - \int
_t^T g\bigl(s,B_{\cdot},
\Yc_s^\P, \Zc_s^\P,
\ah_s\bigr) \dvtx d \langle B \rangle _s
\nonumber
\\[-8pt]
\\[-8pt]
\nonumber
&&{}- \int
_t^T \Zc_s^\P\cdot
\,dB_s - \Nc_T^\P+ \Nc_t^\P,
\end{eqnarray}
whose solution is a triple of $\overline{\F^+}^\P$-progressive
processes, denoted by $(\Yc^{\P}, \Zc^{\P},\break  \Nc^{\P})$, such that
$\Nc
^{\P}$ is a $\overline{\F^+}^\P$-martingale orthogonal to $B$ and
\eqref
{eq:BSDE} holds true $\P$-a.s. We shall assume sufficient conditions
(see Assumption~\ref{assump.gen} below) to guarantee the existence and
uniqueness of the solution to \eqref{eq:BSDE} under every $\P\in\Pc
_W$. In particular, whenever $\P\in\Pc_S$, \eqref{eq:BSDE} turns out
to be a classical BSDE whose solution satisfies $\Nc^{\P} = 0$ and
$\Yc
^{\P}, \Zc^{\P}$ are $\Fbb^{\P}$-progressive. This is due to the fact
that by Lemma~$8.2$ in \cite{stz3}, every probability measures in $\Pc
_S$ satisfies the predictable martingale representation property and
the Blumenthal 0--1 law. This also implies in this case that $\Yc^{\P
}_0$ is a deterministic constant.

The main purpose of the paper is to study the weak
approximation of the following optimization problem:
%
\begin{equation}
\label{eq:supBSDE} Y_0 := \sup_{\P\in\Pc_S}
\Yc^{\P}_0.
\end{equation}

\begin{Remark}
The above problem $Y_0$ in \eqref{eq:supBSDE} is related to the
solution of the following 2BSDE, in the sense that $Y_0$ is the initial
value of the $Y$ component of its solution
\begin{eqnarray*}
Y_t &=& \xi(B_{\cdot}) - \int_t^T
\bigl( g(s,B_{\cdot}, Y_s, Z_s,
\ah_s) \dvtx\widehat a_s \bigr) \,ds\\
&&{} - \int
_t^T Z_s \cdot \,dB_s +
K_T - K_t, \qquad \Pc_S\mbox{-q.s.},
\end{eqnarray*}
which has been introduced by Soner, Touzi and Zhang \cite{Soner_Touzi_Zhang1}.
We also refer to their Section~3.3 for more details,
and simply emphasize here that given the boundedness assumptions we
make below,
it is not necessary in our setting to work on the subset $\mathcal
P^\kappa_H$ of $\mathcal P_S$ introduced in \cite{Soner_Touzi_Zhang1}.
We would also like to comment on the fact that in \cite
{Soner_Touzi_Zhang1}, the solution $(Y,Z)$ is $\F^+$-progressive, while
we defined the solution to the BSDE \reff{eq:BSDE} to be $\overline
{\F
^+}^\P$-progressive. However, thanks to Lemma $2.4$ of~\cite{stz3}, for
any $\P\in\Pc_W$, any $\overline{\F^+}^\P$-progressive process
$X$ has
a $\P$-version $\widetilde X$ which is $\F^+$-progressive, so that this
is not a real difference.
\end{Remark}

We shall impose the following assumptions on the terminal function
$\xi$ and generator function $g$ throughout the paper.
For ease of notation, and since this function will be the main focus
of our paper, we define the function $f\dvtx[0,T] \x\Om\x\R\x\R
^d \x
\S_d^+ \to\R$
\[
f(t,\xb,y,z,u) := g(t,\xb,y,z,u) \dvtx u.
\]

\begin{Assumption}\label{assump.gen}
(i) $\xi\dvtx\Om\longrightarrow\R$ is a bounded Lipschitz
continuous function.\vspace*{-6pt}
\begin{longlist}[(iii)]
\item[(ii)] The process $t \longmapsto f(t,X_{\cdot}, Y_t, Z_t, \nu
_t)$ is progressively measurable given progressive processes $(X, Y, Z,
\nu)$,
and is uniformly continuous with modulus $\rho$ in the sense that for
every $s \le t$ and $\xb, y, z, u$,
\[
\bigl|f(t,\xb_{s \wedge\cdot}, y, z, u)
 - f(s, \xb_{s \wedge\cdot}, y, z, u) \bigr| \le
\rho(t-s).
\]
\item[(iii)] $f$ is uniformly Lipschitz in $(\xb, y,z)$, that is, for
all $(t,\xb_1, \xb_2, y_1,y_2, z_1\break ,z_2,u)$,
\[
\bigl|f(t,\xb_1,y_1,z_1,u)-f(t,
\xb_2,y_2,z_2,u)\bigr| \leq \mu \bigl(\|
\xb_1 - \xb_2\|_t + |y_1-y_2|+|z_1-z_2|
\bigr),
\]
for some constant $\mu>0$.

\item[(iv)] The map $u \longmapsto f(t,\xb,y,z,u)$ is convex and
uniformly continuous for every
$(t,\xb,y,z) \in[0,T] \x\Om\x\R\x\R^d$.

\item[(v)] We have the following integrability condition, for some
constant $C > 0$:
\[
\sup_{(t, \xb, u) \in[0,T] \x\Om\x A} \bigl|f(t, \xb ,0,0,u)\bigr|\leq C.
\]
\end{longlist}
\end{Assumption}

Let us give an existence and equivalence result on the above
2BSDE, whose proof is postponed to Section~\ref{sec:equiv}.

\begin{Theorem} \label{theo:equivalencePsPw}
Suppose that Assumption~\ref{assump.gen} holds true.
Then for every $\P\in\Pc_W$, the BSDE \eqref{eq:BSDE} has a unique
solution $(\Yc^{\P}, \Zc^{\P}, \Nc^{\P})$.
Moreover, we have
%
\begin{equation}
\label{eq:equivPsPw} Y_0 :=\sup_{\P\in\Pc_S}
\Yc^{\P}_0 = \sup_{\P\in\Pc_W} \E^\P
\bigl[ \Yc^{\P}_0 \bigr].
\end{equation}
\end{Theorem}

\begin{Remark} \label{rem:2BSDE_PDE}
Suppose that $\xi(\xb) = \xi_0(\xb_T)$ and $f(t,\xb,y, z,u) =
f_0(t, \xb
_t,\break  y, z, u) $ for some deterministic functions $\xi_0 \dvtx\R^d
\longrightarrow\R$ and $f_0 \dvtx[0,T] \x\R^d \x\R\x\R^d \x
A\longrightarrow\R$.
In this Markovian case, the value function can be given as the
viscosity solution $v(t,x)$ of the nonlinear equation
%
\begin{equation}
\label{eq:PDE} - \partial_t v - \sup_{a \in A} \biggl(
\frac{1}{2} a\dvtx D^2 v - f_0(t, x, v, Dv, a)
\biggr) = 0,
\end{equation}
with terminal condition $v(T,x) = \xi_0(x)$. We refer the reader to
the paper by Soner, Touzi and Zhang \cite{Soner_Touzi_Zhang1} for more
information.
\end{Remark}

\subsection{Weak approximation of 2BSDEs}

Under every probability measure $\P\in\Pc_S$, the canonical process
$B$ is a continuous martingale,
which drives the BSDE~\eqref{eq:BSDE}.
When this martingale is approximated ``weakly'' by a sequence of martingales,
it follows by the robustness property for BSDEs proved by Briand,
Delyon and M\'emin \cite{BriandDelyonMemin}
that the corresponding solutions of the BSDEs driven by the
approximating martingales converge to $\Yc^{\P}$ (see their Theorem $12$).
In the context of 2BSDEs \eqref{eq:supBSDE}, the solution is given as
the supremum of the family of solutions to
BSDEs driven by the family of martingales $( B|_{\P})_{\P\in\Pc_S}$.
Therefore, it is natural, in order to obtain weak approximation
properties, to consider a sequence of families of BSDEs driven by
approximating martingales.
In particular, we shall consider a family of discrete time martingales,
motivated by its application in the numerical approximation described
in Section~\ref{sec:num}.

For every $n \ge1$, we denote by $\Delta_n = (t^n_k)_{0 \le k \le
n}$ a discretization of $[0,T]$, such that $0 = t_0^n < t_1^n < \cdots
< t_n^n = T$.
Let $|\Delta_n| := \sup_{1 \le k \le n} (t_k^n - t^n_{k-1})$, and we
suppose that $|\Delta_n| \longrightarrow0$ as $n \longrightarrow
\infty$.
For ease of presentation, we shall simplify the notation of the time
step size $\Delta t^n_k := t^n_k - t^n_{k-1}$ into $\Delta t$ when
there is no ambiguity. Similarly, we suppress the dependence in $n$ of
$t^n_k$ and write instead $t_k$.

For every $n \ge1$, let $(\Om^n, \Fc^n, \P^n)$ be a probability
space containing $n$ independent random variables $(U_k)_{1 \le k \le n}$.
Moreover, we consider a family of functions $(H^n_k)_{1 \le k \le n, n
\ge1}$ such that every $H_k^n \dvtx A \x[0,1] \longrightarrow\R^d$ is
continuous in $a$ and for some $\delta> 0$, we have for any $a$
%
\begin{eqnarray}
\label{eq:Cond_H} \E \bigl[ H_k^n(a, U_k)
\bigr]& =& 0,\qquad \operatorname{Var}\bigl(H_k^n(a,U_k)
\bigr) = a \Delta t,
\nonumber
\\[-8pt]
\\[-8pt]
\nonumber
 \E \bigl[\bigl |H_k^n(a,U_k)\bigr|^{2+\delta}
\bigr] &\leq& C \Delta t^{1 + {\delta}/{2}},
\end{eqnarray}
where it is understood that the expectation is taken under $\P^n$.

Define the filtration $\F^n:=(\mathcal F^n_{t_k})_{1\leq
k\leq n}$, with $\Fc^n_{t_k} := \sigma(U_1, \ldots, U_k)$ and denote by
$E_n$ the collection of all $\F^n$-predictable $A$-valued processes $e
= (a^e_{t_1}, \ldots,\break  a^e_{t_n})$. Then for every $e \in E_n$, $M^{e}$
is defined by
%
\begin{equation}
\label{eq:M_e_by_H} M^{e}_{t_k} := \sum
_{i \le k} H^n_i\bigl(a^e_{t_i},
U_i\bigr).
\end{equation}

\begin{Remark}
An easy example is when $U_k$ is a Gaussian random vector
($d$-dimension) with distribution $N(0,I_d)$ and $H^n_k(a,u) := a u
\Delta t$.
More examples which induce several different numerical schemes will be
given later in Section~\ref{sec:num}.
\end{Remark}

By abuse of notation, we define a continuous time filtration $\F^n
= \break (\Fc^n_t)_{0 \le t \le T}$,
with $\Fc^n_t := \Fc^n_{t_k}, \forall t \in[t_k, t_{k+1})$
and a continuous time martingales $M^e_t := M^e_{t_k}$, for all $t \in
[t_k, t_{k+1})$ on $(\Om^n, \Fc^n, \P^n)$.
We next consider the completed filtration under~$\P^n$, denoted by
$\mathbb G^n:=\overline{\F^n}^{\P^n}$.
Clearly, $\mathbb G^n$ is right-continuous and complete under $\P^n$,
and $M^e$ is a right-continuous, piecewise constant in time, $\mathbb
G^n$-martingale for every $e \in E_n$.
We notice that the predictable quadratic variation of $M^e$ is given by
\[
\bigl\langle M^e \bigr\rangle _{t_k} = \sum
_{i \le k} \Delta\bigl\langle M^e \bigr\rangle
_{t_k} = \sum_{i
\le k} a^e_i
\Delta t_i.
\]

For every $n \ge1$, with the time discretization $\Delta_n$, we
introduce the truncated generator $f_n(t, \xb, y, z, a) := g_n(t, \xb,
y ,z ,a) \dvtx a$
where
\[
g_n(t, \xb, y, z, a) := g(t_k, \xb, y, z, a)\qquad
\mbox{whenever } t \in[t_k, t_{k+1}).
\]
Then for every $e \in E_n$ and $n \ge1$, we consider the following BSDE:
%
\begin{eqnarray}
\label{eq:Discret_BSDE} \Yc^e_t &=& \xi\bigl(\Mh^e_{\cdot}
\bigr)
 - \int_t^T g_n\bigl(s,
\Mh^e_{\cdot}, \Yc^e_{s^-},
\Zc^e_s, a^e_s\bigr)\dvtx d
\bigl\langle M^e \bigr\rangle _s
\nonumber
\\[-8pt]
\\[-8pt]
\nonumber
&&{} - \int
_t^T \Zc^e_s \cdot
\,dM^e_s -\mathcal N^e_T +
\mathcal N^e_t,
\end{eqnarray}
whose solution is a triple of $\G^n$-progressive processes $(\Yc^e,
\Zc
^e, \mathcal N^e)$
such that $\mathcal N^e$ is a $\G^n$-martingale orthogonal to $M^e$,
and where $\Mh^e$ denotes the continuous interpolation of $M^e$ on the
interval $[0,T]$.
We then have the following wellposedness result for the BSDE \reff
{eq:Discret_BSDE},
which is a direct consequence of Proposition~\ref{prop:sol_GBSDE}
reported in Section~\ref{sec:equiv} and the fact that by taking
conditional expectation with respect to $\G^n$, the component the
solution to \eqref{eq:Discret_BSDE} is given explicitly by the
following scheme:
%
\begin{equation}
\label{eq:Yc} %
\cases{ \Yc^e_{t_n} = \xi\bigl(
\Mh^e_{\cdot}\bigr),\vspace*{2pt}
\cr
\Yc^e_{t_k}
= \E^n_{t_k} \bigl[\Yc^e_{t_{k+1}}\bigr]
- f \bigl(t_k, \Mh ^e_{\cdot
},
\Yc^e_{t_k}, \Zc^e_{t_{k}},
a^e_{t_k} \bigr) \Delta t,\vspace *{2pt}
\cr
\displaystyle \Zc^e_{t_{k}} = \E^n_{t_k} \biggl[
\frac{\Delta\Yc^e_{t_{k+1}}
(a^e_{t_k})^{-1}\Delta M_{k+1}^e}{ \Delta t} \biggr],\vspace*{2pt}
\cr
\Delta\mathcal N^e_{t_{k+1}}=
\Yc^e_{t_{k+1}}-\E^n_{t_k} \bigl[\Yc
^e_{t_{k+1}}\bigr]-\Zc^e_{t_k} \cdot
\Delta M^e_{t_{k+1}},} %
\end{equation}
where $\E^n_{t_k}[\cdot]$ represents the conditional expectation
w.r.t. $\Fc^n_{t_k}$.

\begin{Lemma} \label{lemm:BSDE_e}
Suppose that Assumption~\ref{assump.gen} holds true.
Then for every $n \ge1$ and $e \in E_n$, there is a unique solution
$(\Yc^e, \Zc^e, \Nc^e)$ to the BSDE \eqref{eq:Discret_BSDE} such that
\[
\E^{\P^e} \biggl[ \sup_{0 \le t \le T} \biggl[ \bigl|
\Yc_t^e\bigr|^2 + \int_0^t
\bigl|\bigl(a^e_s\bigr)^{1/2}\Zc^e_s\bigr|^2
\,ds + \bigl\langle\mathcal N^e_t\bigr\rangle \biggr]
\biggr] \le C,
\]
for some constant $C$ independent of $e$ and $n$.
In particular, $\Yc^e_0$ is a deterministic constant.
\end{Lemma}

\begin{pf}
The existence and uniqueness is immediate by \reff{eq:Yc}.
Moreover, Proposition~\ref{prop:sol_GBSDE} gives us the required
estimate for $n\geq n_0$ for some $n_0$. Since only a finite number of
values for $n$ remains, the result is immediate by the fact that the
solution given in \reff{eq:Yc} has the required integrability.
\end{pf}

For every $n \ge1$, denote now
%
\begin{equation}
\label{eq:BSDE_Dn} Y_0^n := \sup_{e \in E_n}
\Yc^e_0.
\end{equation}

The next assumption is a monotonicity condition for the
discretized BSDEs.

\begin{Assumption} \label{assum:monotonicity}
For every $e \in E_n$ and $n \ge1$, the backward scheme in \eqref
{eq:Yc} is monotone,
that is, let $(\Yc^1, \Zc^1)$, $(\Yc^2, \Zc^2)$ be two solutions of
\eqref
{eq:Yc}, then
\[
\Yc^1_{t_{k+1}} \le\Yc^2_{t_{k+1}} \quad
\Longrightarrow\quad \Yc^1_{t_k} \le\Yc^2_{t_k}\qquad
\forall k = 0, \ldots, n-1.
\]
\end{Assumption}

We now state our main result.

\begin{Theorem}\label{theo:weak_convergence}
\textup{(i)} Suppose that Assumption~\ref{assump.gen} holds true. Then
\[
\liminf_{n \to\infty} Y^n_0 \ge
Y_0.
\]
\textup{(ii)} Suppose in addition that Assumption~\ref{assum:monotonicity}
holds and $f$ does not depend on $z$. Then
\[
\lim_{n \to\infty} Y^n_0 = Y_0.
\]
\end{Theorem}

\begin{Remark}
We are not able to show (ii) when the generator depends on $z$. This
is deeply linked to the fact that there are considerable difficulties
to obtain any convergence results for the $z$ part of the solution.
Moreover, since we are working under many measures, the canonical
process is no longer always a Brownian motion, which prevents us from
recovering the strong regularity results of \cite{Zhang2004}, for
instance. We leave this open problem for future research.
\end{Remark}

In the case where $f = 0$, the solution of the 2BSDE is the so
called $G$-expectation of Peng.
Then, in particular, the above result generalizes the weak convergence
result for $G$-expectation in Dolinsky, Nutz and Soner \cite
{Dolinsky_Nutz_Soner}.
We shall report its proof later in Section~\ref{sec:cvg}.

\begin{Remark} \label{rem:monotonicity}
Let $(\Yc^1, \Zc^1)$, $(\Yc^2, \Zc^2)$ be two solutions of \eqref
{eq:Yc}, we have then clearly
%
\begin{eqnarray}
\label{eq:monotonie} &&(1 - L_{t_k,y} \Delta t ) \bigl( \Yc^1_{t_k}
- \Yc^2_{t_k} \bigr)
\nonumber
\\[-8pt]
\\[-8pt]
\nonumber
&&\qquad = \E^n_{t_k}
\bigl[ \bigl( \Yc^1_{t_{k+1}} - \Yc^2_{t_{k+1}}
\bigr) \bigl( 1 + L_{t_k,z} \cdot(\a_{t_k})^{-1}
\Delta M^e_{t_{k+1}} \bigr) \bigr],
\end{eqnarray}
where $L_{t_k,y}$ (resp., $L_{t_k,z}$) is a $\R$-valued
(resp., $\R^d$-valued) and $\Fc^n_{t_k}$-measurable random
variable bounded by the Lipschitz constant $L_{f,y}$ (resp.,
$L_{f,z}$).
Then for $\Delta t$ small enough, the monotonicity condition in
Assumption~\ref{assum:monotonicity} holds whenever
\[
\bigl| L_{f,z} H^n_k(a_{t_k},
U_k) \bigr| \le|a_{t_k}|\qquad \forall1 \le k \le n.
\]
In particular, when $f$ is independent of $z$, Assumption~\ref
{assum:monotonicity} always holds true for $\Delta t$ small enough.
\end{Remark}

\section{Numerical schemes for 2BSDEs}
\label{sec:num}

As discussed in Remark~\ref{rem:2BSDE_PDE}, the solution of the
Markovian 2BSDE \eqref{eq:supBSDE} can be given as viscosity solution
of a parabolic fully nonlinear PDE, for which a comparison principle holds.
Several monotone numerical schemes have been proposed for PDEs in or
closed to this form,
for example, the generalized finite difference scheme of Bonnans,
Ottenwaelter and Zidani \cite{Bonnans_2004},
the semi-Lagrangian scheme of Debrabant and Jakobsen~\cite{Debrabant_Jakobsen},
and the probabilistic scheme of Fahim, Touzi and Warin \cite
{Fahim_Touzi_Warin}, Guo, Zhang and Zhuo \cite{GuoZhangZhuo}, where the
convergence is ensured by the monotone convergence theorem of Barles
and Souganidis \cite{Barles_Souganidis}.

Similar to Tan \cite{TanNumStoContr} in the context of
non-Markovian control problems,
we can interpret these schemes as a system of controlled Markov chains.
Using these controlled Markov chains as the families of driving
martingale $(M^e)_{e \in E_n}$ in \eqref{eq:Discret_BSDE}, Theorem~\ref
{theo:weak_convergence} also justifies the convergence of the
corresponding numerical schemes.
Moreover, it permits to extend these numerical schemes to the
non-Markovian case.
The aim of this section is to present a general abstract numerical
scheme for 2BSDEs, which we then specialize in two particular examples.
In particular, these schemes are coherent with the numerical methods
proposed and tested in the previous literature, for which we can refer
to \cite{Fahim_Touzi_Warin,TanSplitting,GuoZhangZhuo}, etc.
We nonetheless start by studying the solution to the discrete-time BSDEs.

\subsection{An explicit scheme}

We notice that for every fixed $e \in E_n$ and $n \ge1$,
the backward iteration in \eqref{eq:Yc} is in fact the so called
implicit scheme for BSDEs.
In practice, we consider also the following explicit scheme:
%
\begin{equation}
\label{eq:Yct} %
\cases{ \Yct^e_{t_n} = \xi\bigl(
\Mh^e_{\cdot}\bigr),\vspace*{2pt}
\cr
\Yct^e_{t_k}
= \E^n_{t_k} \bigl[\Yct^e_{t_{k+1}}\bigr]
- f \bigl(t_k, \Mh ^e_{\cdot
},
\E^n_{t_k}\bigl[\Yct^e_{t_{k+1}}\bigr],
\Zct^e_{t_{k+1}}, a^e_{t_k} \bigr) \Delta
t,\vspace*{2pt}
\cr
\displaystyle\Zct^e_{t_{k}} = \E^n_{t_k}
\biggl[ \frac{\Delta\Yc^e_{t_{k+1}}
(a^e_{t_k})^{-1}\Delta M_{k+1}^e}{ \Delta t} \biggr].} %
\end{equation}
Denote
%
\begin{equation}
\label{eq:Ytn} \Yt^n_0 := \sup_{e \in E_n}
\Yct^e_0.
\end{equation}

The following lemma shows that the implicit and explicit schemes
only differ by an amount proportional to $\Delta_n$.

\begin{Lemma}
There is a constant $C$ independent of $n \ge1$ such that
\[
\bigl|Y^n_0 - \Yt^n_0\bigr|\le C |
\Delta_n|.
\]
\end{Lemma}

\begin{pf} It is enough to prove that there is some constant $C > 0$
independent of $n \ge1$ and $e \in E_n$ such that
\[
\bigl|\Yc^e_0 - \Yct^e_0\bigr| \le C |
\Delta_n|.
\]
First, by \eqref{eq:Yc} and \eqref{eq:Yct} and the Lipschitz property
of the generator $f$, it is clear that for every $0 \le k \le n-1$,
there are bounded $\mathcal G^n_{t_k}$-random variables $\alpha_k$ and
$\beta_k$ such that
\begin{eqnarray*}
\bigl( \Yc^e_{t_k} - \Yct^e_{t_k}
\bigr) &=& \E^n_{t_k} \bigl[ \Yc^e_{t_{k+1}}
- \Yct ^e_{t_{k+1}} \bigr] + \alpha_k \bigl(
\Yc^e_k - \E^n_{t_k}\bigl[
\Yct^e_{t_{k+1}}\bigr]\bigr) + \beta_k \cdot\bigl(
\Zc^e_{t_{k+1}} - \Zct^e_{t_{k+1}}\bigr)
\\
&=& (1+ \alpha_k \Delta t) \E^n_{t_k} \bigl[
\bigl( \Yc^e_{t_{k+1}} - \Yct^e_{t_{k+1}}
\bigr) \bigl( 1 + (1 + \alpha_k \Delta t)^{-1}
\beta_k \cdot\Delta M_{k+1} \bigr) \bigr]
\\
&&{} + f\bigl(t_k, \Mh^e_{\cdot},
\Yc^e_{t_k}, \Zc^e_{t_{k+1}},
a^e_{t_k}\bigr) \Delta t^2.
\end{eqnarray*}
Then using the Young inequality $(a+b)^2 \le(1 + \gamma h) a^2 + (1 +
\frac{1}{\gamma h}) b^2$ and the Cauchy--Schwarz inequality, we get for
some constant $C$ independent of $e$ and $k$,
\begin{eqnarray*}
\bigl( \Yc^e_{t_k} - \Yct^e_{t_k}
\bigr)^2 &\le& (1+\gamma\Delta t) (1 + C \Delta t)
\E^n_{t_k} \bigl[ \bigl(\Yc^e_{t_{k+1}} -
\Yct^e_{t_{k+1}}\bigr)^2 \bigr]
\\
&&{}+ C f^2\bigl(t_k, \Mh^e_{\cdot},
\Yc^e_{t_k}, \Zc^e_{t_{k+1}},
a^e_{t_k}\bigr) \Delta t^2.
\end{eqnarray*}
Taking expectations en each side and using the Lipschitz property of
$f$, we get
\begin{eqnarray*}
\E^n_0 \bigl[ \bigl( \Yc^e_{t_k} -
\Yct^e_{t_k} \bigr)^2 \bigr]& \le &
\E^n_0 \bigl[ (1 + C \Delta t) \bigl(\Yc^e_{t_{k+1}}
- \Yct ^e_{t_{k+1}} \bigr)^2 \bigr]
\\
&&{} + C\Delta t^2 \E^n_0 \bigl[\bigl|
\Mh^e\bigr|^2 + \bigl|\Yc^e\bigr|^2 + \bigl|\Zc
^e\bigr|^2 \bigr].
\end{eqnarray*}
Finally, it is enough to conclude using the Gronwall lemma together with
the estimates given by Lemma~\ref{lemm:BSDE_e}.
\end{pf}

For every $n \ge1$, we can reformulate the problem \eqref
{eq:BSDE_Dn} for $Y^n_0$ and \eqref{eq:Ytn} for $\Yt^n_0$ as a
numerical scheme defined on
\[
\Lambda^n := \bigcup_{0 \le k \le n} \{
t_k \} \x\R^{d \x(k+1)}.
\]
For every $n \ge1$, $(t_k, \xb) \in\Lambda^n$ and $a \in A$, we
define $M^{t_k, \xb, a} \in\R^{d \x(k+2)}$ by
\[
\cases{ M^{t_k, \xb, a}_{t_i} := \xb_i,
& \quad$\mbox{for every $i \le k$},$\vspace *{2pt}
\cr
M^{t_k, \xb, a}_{t_{k+1}} :=
M^{t_k, \xb, a}_{t_k} + H_{k+1}^n(a,
U_{k+1}).} %
\]
We then define $u^n \dvtx\Lambda_n\longrightarrow\R$ and $\ut^n
\dvtx\Lambda
_n \longrightarrow\R$ by the following backward iterations.
The terminal conditions are given by
\[
u^n(t_n, \xb) := \ut^n(t_n, \xb)
:= \xi(\hat{\xb})\qquad \forall \xb\in\R^{d \x(n+1)},
\]
and the backward iteration for $u^n$ and $\tilde u^n$ are given by,
for all $\xb\in\R^{d \x(k+1)}$,
%
\begin{equation}\qquad
\label{eq:uan} %
\cases{ \displaystyle u^n(t_k, \xb) = \sup
_{a \in A} u_a^n(t_k, \xb),
\vspace*{2pt}
\cr
u_a^n(t_k, \xb) = \E\bigl[ u
\bigl(t_{k+1}, M^{t_k,\xb,a} \bigr) \bigr] - f\bigl(t_k,
\hat\xb, u_a^n(t_k, \xb),
Du_a^n(t_k, \xb), a \bigr) \Delta t, \vspace
*{2pt}
\cr
\displaystyle Du_a^n(t_k, \xb) := \E \biggl[
\frac{ u(t_{k+1}, M^{t_k,\xb,a} )
a^{-1}\Delta M^{t_k, \xb,a}_{k+1}}{ \Delta t} \biggr],} %
\end{equation}
and
%
\begin{equation}\qquad
\label{eq:abstract_num_scheme} %
\cases{ \displaystyle\ut^n(t_k, \xb) =
\sup_{a \in A} \bigl( \ut_a^n(t_k,
\xb) - f\bigl(t_k, \hat\xb, \ut_a^n(t_k,
\xb), D \ut _a^n(t_k, \xb), a \bigr) \Delta
t \bigr),\vspace*{2pt}
\cr
\ut_a^n(t_k, \xb)
:= \E\bigl[ \ut\bigl(t_{k+1}, M^{t_k,\xb,a} \bigr) \bigr],\vspace
*{2pt}
\cr
\displaystyle D \ut_a^n(t_k, \xb) := \E \biggl[
\frac{ \ut(t_{k+1}, M^{t_k,\xb
,a} )
a^{-1}\Delta M^{t_k, \xb,a}_{k+1}}{ \Delta t} \biggr].} %
\end{equation}

We have the following dynamic programming result.

\begin{Proposition}\label{prop:truc}
Let Assumption~\ref{assum:monotonicity} hold true, then
\[
\ut^n(0,\0) = \Yt^n_0 \quad\mbox{and}\quad
u^n(0, \0) = Y_0^n.
\]
\end{Proposition}

\begin{pf}It is in fact a standard result from the dynamic programming
principle; see, for example, Bertsekas and Shreve \cite
{BertsekasShreve} for a detailed presentation on this subject. We also
notice that the arguments are almost the same in Theorem~3.4 of Tan
\cite{TanNumStoContr} for a similar problem.
\end{pf}

\subsection{Concrete numerical schemes of 2BSDE}
\label{subsec:num_scheme}

By constructing the driving martingales $(M^e)_{e \in E_n}$ as a
family of controlled Markov chain, we can also compute the solution of
\eqref{eq:BSDE_Dn} using a backward iteration, under some monotonicity
conditions. In particular, it can be considered as a numerical scheme
for the 2BSDE~\eqref{eq:supBSDE}. For particular choices of functions
$(H_k^n)_{1 \le k \le n, n \ge1}$, we may obtain some numerical
schemes, including a finite difference scheme and a probabilistic scheme.

\subsubsection{Finite difference scheme}

Let us stay in the one-dimensional case $d=1$ for notational
simplicity, where $\Delta x \in\R$ is the parameter of the space
discretization. Denote $p_a := a \Delta t /\Delta x^2$, suppose that $
p_a\le1/2$ for all $a \in A$.
Clearly, for every $n \ge1$ and space discretization $\Delta x$, we
can construct a function $H^n\dvtx A \x[0,1] \longrightarrow\{ -
\Delta
x, 0, \Delta x \} $ such that, for any uniformly distributed random
variable~$U$
\[
\P^n \bigl[ H^n(a,U) = \Delta x \bigr] =
\P^n \bigl[ H^n(a,U) = -\Delta x \bigr] = p_a,
\]
and
$\P^n  [H^n(a,U) = 0  ] = 1 - 2p_a$.
Let $H_k^n := H^n$, and denote $\xb^{k, \pm} := (x_0, \ldots, \break  x_k, x_k
\pm\Delta x)$ and $\xb^{k,0} = (x_0, \ldots, x_k, x_k)$ for every
$\xb
= (x_0, \ldots, x_k)$.
Then it follows by a direct computation that the numerical iteration
in \eqref{eq:abstract_num_scheme} turns to be
%
\begin{equation}
\label{eq:FD_scheme} \ut(t_k, \xb) := \ut\bigl(t_{k+1},
\xb^{k,0}\bigr) + \sup_{a \in A} \biggl\{
\frac
{1}{2} a D^2 \ut- f (\cdot, \ut_a, D \ut, a )
(t_k, \xb) \biggr\},
\end{equation}
where $ \ut_a(t_k, \xb) = \ut(t_{k+1}, \xb^{k,0}) + \frac{1}{2} a
\Delta t D^2 \ut(t_k, \xb)$, with
\[
D^2 \ut(t_k, \xb) = \frac{\ut(t_{k+1}, \xb^{k,+}) - 2 \ut
(t_{k+1}, \xb
) + \ut(t_{k+1}, \xb^{k,-})}{\Delta x^2}
\]
and
\[
D \ut(t_k, \xb) = \frac{\ut(t_{k+1}, \xb^{k,+}) - \ut(t_{k+1}, \xb^{k,-})}{2 \Delta x}.
\]

\begin{Remark}
(i) For the above choice of $(H^n_k)_{1 \le k \le n}$, Assumption~\ref
{assum:monotonicity} holds true whenever $\Delta x \le L_{f,z}$.\vspace*{-6pt}
\begin{longlist}[(iii)]
\item[(ii)]To ensure that $p_a := a \Delta t /\Delta x^2 \le1/2$, we
should choose $\Delta x \sim\sqrt{\Delta t}$. Moreover, the family of
functions $(H_k^n)_{0 \le k \le n}$ associated with the finite
difference scheme satisfies condition \eqref{eq:Cond_H}.

\item[(iii)] In the high dimensional case $d > 1$, the construction of
finite difference scheme will be harder in general. We refer to Kushner
and Dupuis \cite{KushnerDupuis} in the case where all $a \in A$ are
diagonal dominant, and also to Bonnans, Ottenwaelter and Zidani
\cite{Bonnans_2004} in general cases.
\end{longlist}
\end{Remark}

\subsubsection{Probabilistic scheme}

For parabolic nonlinear PDEs including \eqref{eq:PDE}, Fahim, Touzi
and Warin \cite{Fahim_Touzi_Warin} proposed a probabilistic scheme,
which was reinterpreted and generalized in a non-Markovian stochastic
control context in Tan~\cite{TanNumStoContr}. We can easily adapt this
probabilistic scheme in our context.

Let $a_0 \in\S_d^+$ be a fixed constant, denote $\sigma_0 =
a_0^{1/2}$. Suppose that for all $a \in A$,
\[
a \ge a_0 \quad\mbox{and}\quad 1 - \tfrac{1}{2} (a - a_0)
a_0^{-1} \ge0.
\]
For every $n \ge1$, denote $\rho_n \dvtx A \x\R^d \longrightarrow
\R$ by
%
\begin{equation}
\label{eq:def_f} \rho_n( a, x) := \frac{1}{ (2 \pi\Delta t)^{d/2} | \sigma_0 |^{1/2}} \exp \biggl( -
\frac
{1}{2} \Delta x^{-1} x^T a_0^{-1}
x \biggr) \eta_n(a,x),
\end{equation}
with
\[
\eta_n(a,x) := \bigl( 1- \tfrac{1}{2} a \cdot
a_0^{-1} + \tfrac{1}{2} \Delta t^{-1} a
\cdot a_0^{-1} x x^T \bigl(a_0^T
\bigr)^{-1} \bigr).
\]
It is easy to verify that $x \longmapsto\rho_n(a,x)$ is a probability
density function for every $a \in A$. Then following Tan \cite
{TanNumStoContr}, we can construct $H^n(a,x)$ which is continuous in
$a$ and such that $H^n(a,U)$ is a random variable of density function
$\rho_n(a,x)$ whenever $U \sim\Uc[0,1]$.

To make Assumption~\ref{assum:monotonicity} hold true, we suppose
in addition that $f$ is independent of $z$ (see Remark~\ref{rem:monotonicity}).
Define then the family of functions $(H^n_k)_{ 1 \le k \le n}$ by
$H^n_k = H^n$.
We can then rewrite $\ut_a^n$ in \eqref{eq:abstract_num_scheme} by the
following: let $\Delta W \sim N(0, \Delta t I_d)$,
\begin{eqnarray*}
\ut_a^n(t_k, \xb) &=& \E \bigl[ \ut
\bigl(t_{k+1}, \bigl(\xb, x_k + H^n(a, U)\bigr)
\bigr) \bigr]
\\
&=& \E \bigl[ \ut\bigl(t_{k+1}, (\xb, x_k +
a_0 \Delta W) \bigr) \eta_n(a, a_0 \Delta W)
\bigr]
\\
&=& \E \bigl[ \ut \bigl(t_{k+1}, (\xb, x_k +
a_0 \Delta W) \bigr) \bigr]
\\
&&{} + \frac{1}{2} \Delta t a \cdot\E \biggl[ \ut \bigl(t_{k+1}, (
\xb, x_k + a_0 \Delta W) \bigr) \bigl(
\sigma_0^T\bigr)^{-1} \\
&&\hspace*{65pt}{}\times\frac{\Delta W_{k+1}
\Delta
W^T_{k+1} - \Delta t I_d}{\Delta t^2}
\sigma_0^{-1} \biggr].
\end{eqnarray*}
Therefore, the explicit numerical scheme \eqref{eq:abstract_num_scheme}
can be rewritten in the following way:
in a probability space $(\Om^0, \Fc^0, \P^0)$, let $X^0 := (a_0
W_{t_0}, \ldots, a_0 W_{t_n}) \in\R^{d \x(n+1)}$, where $W$ is a
standard $d$-dimensional Brownian motion.
Let $\widehat{X}$ denote continuous time process obtained by linear
interpolation of the discrete time process $X^0$.
The terminal condition is given by $\Yt_{t_n} = \xi(\widehat
{X}_{\cdot
})$, and the backward iteration:
%
\begin{equation}
\label{eq:NumScheme2} \Yt_{t_k} := \E_{t_k} [\Yt_{t_{k+1}}] +
\Delta t G \bigl(t_k, \widehat {X^0}_{\cdot},
\E_{t_k} [\Yt_{k+1}], \Gamma_{t_k} \bigr),
\end{equation}
with
\[
\Gamma_{t_k} := \E_{t_k} \biggl[ \Yt_{t_{k+1}} \bigl(
\sigma_0^T\bigr)^{-1} \frac{\Delta W_{k+1}
\Delta W^T_{k+1} - \Delta t I_d}{\Delta t^2}
\sigma_0^{-1} \biggr]
\]
and
\[
G(t, \xb, y, \gamma) := \sup_{a \in U} \biggl( f\biggl(t,\xb, y +
\frac
{1}{2} a \cdot\gamma\Delta t, a\biggr) + \frac{1}{2} a \cdot
\gamma \biggr).
\]
Notice that the above scheme is closely related to the scheme proposed
by Fahim, Touzi and Warin \cite{Fahim_Touzi_Warin} for nonlinear PDEs.

\subsection{Numerical examples}\label{sec:num2}

We provide here some numerical tests on the schemes proposed in
Section~\ref{subsec:num_scheme}.
Let $d = 1$ for simplicity, we shall consider two different equations
with the following generators $f_1$ and $f_2$:
%
\begin{eqnarray}
\label{eq:2BSDE_ex1} f_1(t,x, y, z, a)& := & \inf_{r \in K}
\{r y a \}\qquad \mbox{for some compact set } K \subset\R,
\\
\label{eq:2BSDE_ex2} f_2(t,x,y,z,a) &:= & \tfrac{1}{2} \bigl( (
\sqrt{a}z + b /\sqrt{a} )^- \bigr)^2 - z b - \tfrac{1}{2}
b^2/a,
\end{eqnarray}
and the terminal condition is given by
%
\begin{equation}
\label{eq:2BSDE_tc} \xi( \xb) := K_1 + \biggl(\int_0^T
\xb_t \,dt- K_1 \biggr)^+ - \biggl( \int
_0^T \xb_t \,dt - K_2
\biggr)^+,
\end{equation}
for some constant $K_1 \le K_2$.

We would like to point out to the reader that the first
example of second-order BSDE with generator \eqref{eq:2BSDE_ex1} is
motivated by a differential game type of problem
\[
\sup_{a \in A} \inf_{r \in K} \E \biggl[ \exp
\biggl( \int_0^T r_s
a_s \,ds \biggr) \xi \bigl(X^a_{\cdot} \bigr)
\biggr],
\]
while the second example with generator \eqref{eq:2BSDE_ex2} is taken
from the robust utility maximization problem studied by Matoussi,
Possama\"{i} and Zhou \cite{mpz} (see the generator in their Theorem~4.1,
when the set $A_a$ is chosen to be $[0,+\infty)$).
We also insist on the fact that the generator $f_2$ depends on the $z$
variable and is of quadratic growth, so that our general convergence
result does not apply in this setting.
Nonetheless, as shown by the numerical results below, our numerical
schemes still converge in this case, leading us to the natural
conjecture that convergence also holds in this more general setting.

Moreover, with the above terminal condition \eqref
{eq:2BSDE_tc}, by adding the variable $M_\cdot:=\int_0^\cdot\xb
_t\,dt$ in
the diffusion system, we can also characterize the solution of the
2BSDE by the following degenerate PDE on $v\dvtx(t,,x,m) \in[0,T] \x
\R^2
\longrightarrow\R$:
%
\begin{equation}
\label{eq:PDE_ex} \partial_t + x \partial_m v + \sup
_{a \in A} \biggl( \frac{1}{2} a \partial^2_{xx}
v + f(t,x, v, \partial_x v, a) \biggr) = 0,
\end{equation}
with terminal condition $v(T, x, m) := K_1 + (m -K_1)^+ - (m- K_2)^+$.

For each of the two 2BSDEs, we implemented the finite
difference scheme given by \eqref{eq:FD_scheme} and the probabilistic
scheme \eqref{eq:NumScheme2}.
As a comparison, we also implemented PDE \eqref{eq:PDE_ex} with a
splitting finite difference scheme,
that is to split it into two PDEs:
\[
\partial_t + x \partial_m v = 0 \quad\mbox{and}\quad
\partial_t + \sup_{a \in A} \biggl( \frac{1}{2} a
\partial^2_{xx} v + f(t,x, v, \partial_x v,
a) \biggr) = 0,
\]
and then to solve the two PDEs sequentially with classical
finite-difference scheme.
Since each equation is one-dimensional, the associated classical
finite-difference scheme is bound to be a good benchmark for our schemes.
We implemented the numerical schemes on a computer with 2.4~GHz CPU and
4G memory.

\begin{figure}

\includegraphics{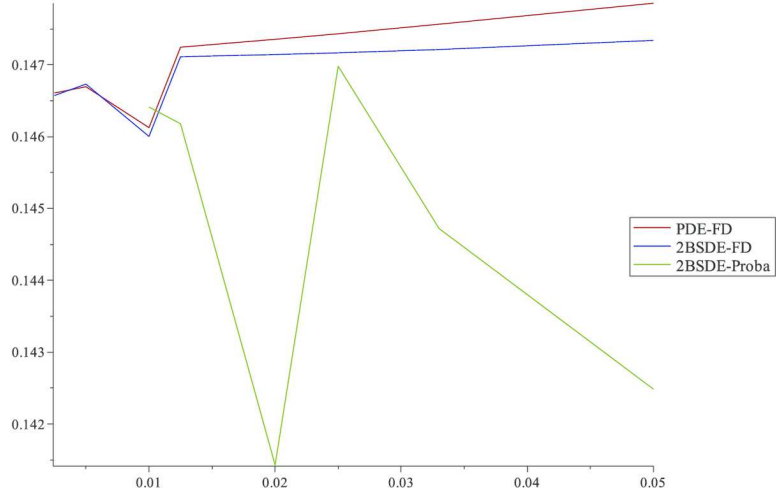}

\caption{The comparison of numerical solutions for 2BSDE with
generator \protect\eqref{eq:2BSDE_ex1}.
The faire value should be very closed to 0.146, and the probabilistic
scheme seems more volatile comparing to the other schemes.}
\label{fig:Asian1}
\end{figure}

In the following two low-dimensional examples, we choose
$X_0 = 0.2$, $K = [-1, 1]$, $K_1 = -0.2$, $K_2 = 0.2$ and $A = [0.04,
0.09]$, corresponding to a volatility uncertainty in $[0.2, 0.3]$.
Using difference time-discretization with time step $\Delta t$,
the numerical solutions of schemes \eqref{eq:FD_scheme} and \eqref
{eq:NumScheme2} are quite stable and closed to the PDE numerical
results w.r.t. the relative error.
In Figures \ref{fig:Asian1} and \ref{fig:Asian2} below, we give the
numerical solutions with different time discretization.
The line PDE-FD denotes the splitting finite-difference method on the
PDE \eqref{eq:PDE_ex},
2BSDE-FD denotes the finite-difference scheme \eqref{eq:FD_scheme} on
the 2BSDE,
and 2BSDE-Proba refers to the probabilistic scheme \eqref
{eq:NumScheme2} on the 2BSDE.
For the probabilistic scheme, we use a simulation-regression to
estimate the conditional expectation arising in the backward iteration
\eqref{eq:NumScheme2}.
When $\Delta t = 0.02$,
a single computation takes 1.72 seconds for PDE-FD, 1.92 seconds for
2BSDE-FD, and 103.2 seconds for the 2BSDE-Proba method (using $2 \x
10^5$ simulations in the simulation-regression method).
In this two-dimensional case, it is not surprising that the
finite-difference scheme is much less time-consuming comparing to the
probabilistic scheme.

\begin{figure}

\includegraphics{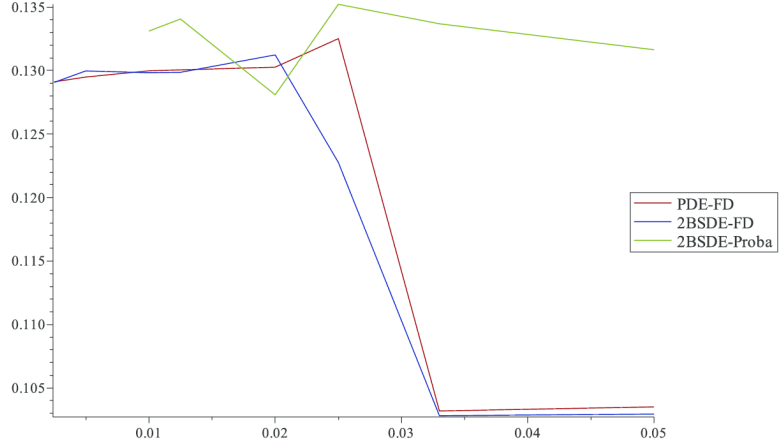}

\caption{The comparison of numerical solutions for 2BSDE with
generator \protect\eqref{eq:2BSDE_ex2}.
The faire value should be closed to 0.129.
For finite-difference scheme, when $\Delta t$ is greater than 0.025,
we need to use a coarser space-discretization to ensure the
monotonicity (similar to the classical CFL condition), which makes a
big difference to the numerical solutions for the case $\Delta t < 0.25$.
However, the convergence as $\Delta t \to0$ is still obvious.}
\label{fig:Asian2}
\end{figure}

\section{Proof of the convergence result}
\subsection{Proof of Theorem \texorpdfstring{\protect\ref{theo:equivalencePsPw}}{2.3}}
\label{sec:equiv}
(i) The wellposedness of the BSDE \reff{eq:BSDE} is a already proved
in Proposition~\ref{prop:sol_GBSDE}.

(ii) We fix a filtered probability space $(\Om^0, \Fc^0,\F^0, \P
^0)$, where the filtration $\F^0$ satisfies the usual hypotheses.\vadjust{\goodbreak} Let
$\Pc_h$ denote the collection of all martingale probability measures
$\P\in\Pc_W$ such that the density process $\hat a$ is piecewise
constant, that is to say $\hat a_t = \sum_{k = 1}^n a_{t_k} 1_{t \in
[t_k, t_{k+1})}$, $d \P\x \,dt$-a.e., for some time discretization $0 =
t_0 < \cdots< t_n = T$.
Let $M$ be a $\F^0$-martingale, whose distribution lies in $\Pc_W$. We
can approximate $M$ by a sequence $\widehat M^n$ such that $\P^0 \circ
(\widehat M^n)^{-1} \in\Pc_h$ and
\[
\mathbb E^\mathbb P \biggl[\int_0^T\bigl|{a^n_t}-{a_t}\bigr|^2
\,dt \biggr]\mathop{\longrightarrow}_{n\rightarrow+\infty} 0\quad\mbox{and} \quad\mathbb E^\mathbb
P \Bigl[\sup_{0\leq t\leq T}\bigl|\widehat{M}^n_t -
M_t\bigr|^2 \Bigr]\mathop{\longrightarrow}_{n\rightarrow+\infty}0,
\]
where $a^n_t := \frac{d \langle\widehat{M}^n \rangle_t}{dt}$ and $a_t
:= \frac{d \langle M \rangle_t}{dt}$.
Then in the spirit of Proposition~\ref{prop:density_approx}, we have
\[
\sup_{\P\in\Pc_h} \E\bigl[ \Yc_0^{\P} \bigr]
= \sup_{\P\in\Pc_W} \E\bigl[ \Yc_0^{\P}
\bigr].
\]
Further, we claim that for every $\P_h \in\Pc_h$,
%
\begin{equation}
\label{eq:Pch} \E\bigl[ \Yc_0^{\P_h}\bigr] \le
Y_0 := \sup_{\P\in\Pc_S} \Yc_0^{\P}.
\end{equation}
It follows that
\[
\sup_{\P\in\Pc_W} \E\bigl[ \Yc_0^{\P} \bigr]
= \sup_{\P\in\Pc_h} \E\bigl[ \Yc_0^{\P}
\bigr] \le \sup_{\P\in\Pc_S} \Yc_0^{\P}.
\]
By the trivial inequality $Y_0 \le\sup_{\P\in\Pc_W} \E[ \Yc
_0^{\P
}]$, we get \eqref{eq:equivPsPw}.

(iii) It remains now to prove the claim \eqref{eq:Pch}. We follow
closely the randomization argument
in Step 3 of the proof of Proposition~3.5 of Dolinsky, Nutz and Soner
\cite{Dolinsky_Nutz_Soner}. We emphasize that the proof in \cite
{Dolinsky_Nutz_Soner} only uses the fact that the set where the density
of the quadratic variation of the canonical process is both convex and
compact, which is the case for our set $A$.
We notice that under $ \P_h \in\Pc_h $, the canonical process $B$ is
a martingale such that the density of its quadratic variation is
piecewise constant. Let us denote it by
\[
\alpha_t :=\sum_{k=0}^{n-1}{
\mathbf1}_{[t_k,t_{k+1})} (t) \alpha(k),
\]
where the $\alpha(k)$ are $\mathcal F_{t_k}$-measurable.
Further, denote
$W_t := \int_0^t \alpha_s^{-1/2} \,dB_s$, which is clearly a $\P
_h$-Brownian motion.
Then by exactly the same arguments as in the step $3$ of the proof of
Proposition~3.5 of \cite{Dolinsky_Nutz_Soner},
we can consider a probability space $(\tilde\Om, \tilde\Fc, \tilde
\P
)$ equipped with a Brownian motion $\tilde W$ and i.i.d. uniformly
distributed r.v. $(\tilde U_k)_{ 1 \le k \le n}$, independent of
$\tilde W$, and construct, using regular conditional probability
measures, random variables $\tilde\alpha(k)$ which are $\sigma(
\tilde
U_j, 1 \le j \le k)\vee\sigma(\tilde W_s, 0\leq s\leq
t_k)$-measurable and such that the following equality holds:
\[
\begin{tabular}{p{250pt}@{}}
the law of $ \bigl(\tilde W, \bigl(\tilde\alpha(i)\bigr)_{0\leq i\leq n-1}
\bigr)$ under $\tilde\P$ $=$ the law of $ \bigl( W, \bigl(\alpha (i)
\bigr)_{0\leq
i\leq n-1} \bigr)$ under $\P_h$.
\end{tabular}
\]
Define next the martingale
\[
\tilde M_t:=\int_0^t \Biggl(\sum
_{k=0}^{n-1}{\mathbf1}_{[t_k,t_{k+1})}(s)
\bigl(\tilde\alpha(k) \bigr)^{1/2} \Biggr)\,d\tilde W_s,\qquad
\tilde\P\mbox{-a.s.}
\]
We deduce that $\tilde\P\circ\tilde M^{-1} = \P_h$. Let us moreover
consider the family of conditional probability measures $(\tilde\P
_c)_{c \in[0,1]^n}$ of $\tilde\P$ w.r.t. the sub-$\sigma$-field
$\sigma( \tilde U_k, 1 \le k \le n)$ and define $\P_c := \tilde\P_c
\circ\tilde M^{-1}$. We have that $\P_c \in\Pc_S$ for every $c \in[0,1]^n$.
It follows that
\[
\E\bigl[ \Yc_0^{\P_h}\bigr] \le \sup_{c \in[0,1]^n}
\Yc_0^{\P_c} \le \sup_{\P\in\Pc_S}
\Yc_0^{\P},
\]
which justifies the claim \eqref{eq:Pch}, and we hence complete the
proof. 

\subsection{Proof of the convergence theorem}
\label{sec:cvg}
To prove Theorem~\ref{theo:weak_convergence}, we shall first provide
some technical lemmas.

\begin{Lemma} \label{lemm:M_e}
Let the functions $H^n_k$ satisfy \eqref{eq:Cond_H}, then there are
some constants $\delta> 0$ and $C>0$ such that for every $e \in E_n$,
$n \ge1$,
%
\begin{equation}
\label{eq:M_e} \E \bigl[\bigl | M^e_{\cdot}\bigr |^{2+\delta}
\bigr] \le C \quad\mbox{and}\quad\bigl |a^e_k\bigr| \le C\qquad \forall1 \le k
\le n, \P^n\mbox{-a.s.}
\end{equation}
In particular this implies that for every $e\in E_n$,
%
\begin{equation}
\label{eq:VQM_e} \bigl\langle M^e\bigr\rangle _t - \bigl
\langle M^e\bigr\rangle _s\leq C\bigl((t-s)+|\Delta
_n|\bigr)I_d,\qquad \mathbb P^n\mbox{-a.s.}
\end{equation}
Moreover, any sequence $(M^{e_n})_{n \ge1}$, with $e_n \in E_n,
\forall n \ge1$, is relatively compact and any limit of the sequence
lies in $\Pc_W$.
\end{Lemma}

\begin{pf} Let $ n \ge1$ and $0 \le s < t \le T$, we can suppose without
loss of generality that $ t -s > |\Delta_n|$ by \eqref{eq:Cond_H}.
Then for every $ e \in E_n$,
\begin{eqnarray*}
&& \E \Bigl[ \sup_{s \le r \le t} \bigl( M_r^e
- M^e_s \bigr)^{2+\delta} \Bigr]
\\
&&\qquad\le C \E \bigl[ \bigl( \bigl[M^e\bigr]_t -
\bigl[M^e\bigr]_s \bigr)^{1 + \delta/2} \bigr]
\\
&&\qquad\le C \E \biggl[ \biggl( \sum_{s \le t_i \le t}
\bigl|H_i^n(e_i, U_i)\bigr|^2
\biggr)^{1 + \delta/2} \biggr]
\\
&&\qquad= C (t-s)^{1 + \delta/2} \E \biggl[ \biggl( \frac{\Delta t}{(t-s)} \sum
_{s \le t_i \le t} \bigl|H_i^n(e_i,
U_i)/\sqrt{\Delta t}\bigr|^2 \biggr)^{1 + \delta/2} \biggr]
\\
&&\qquad\le C (t-s)^{1 + \delta/2} \E \biggl[ \frac{\Delta t}{(t-s)} \sum
_{s \le t_i
\le t} \bigl|H_i^n(e_i,
U_i)/\sqrt{\Delta t}\bigr|^{2 + \delta} \biggr]
\\
&&\qquad\le C (t-s)^{1+\delta/2},
\end{eqnarray*}
where the first inequality follows from BDG inequality, the second
from Jensen's inequality and the last from \eqref{eq:Cond_H}.
It follows that \eqref{eq:M_e} and \eqref{eq:VQM_e} hold true,
and hence any sequence $(M^{e_n})_{n \ge1}$ such that $e_n \in E_n$
is relatively compact (see, e.g., Stroock and Varadhan \cite{Stroock_1979}).
Finally, let $\P$ be a limit probability measure, it follows by
exactly the same argument as in Lemma~3.3 of Dolinsky, Nutz and Soner
\cite{Dolinsky_Nutz_Soner} that $\P\in\Pc_W$, which completes the
proof.
\end{pf}

\begin{Lemma} \label{lemm:Lip_un}
Let $u^n$ be defined in \eqref{eq:abstract_num_scheme} and \eqref
{eq:uan}, then there is constant $C$ independent of $n$ such that
%
\begin{equation}
\label{eq:un_Lip} \bigl| u^n(t_k, \xb_1) -
u^n\bigl(t_{k+1}, \xb_2^{t_k}\bigr) \bigr|
\le C \bigl( |\xb_1 - \xb_2|_k + \sqrt{|
\Delta_n|} \bigr).
\end{equation}
\end{Lemma}

\begin{pf} {(i)} Suppose that $u^n(t_{k+1}, \xb)$ is Lipschitz in $\xb$ with
Lipschitz constant $L_{k+1}$, let $\xb^1, \xb^2 \in\R^{d \x k}$,
then using the same argument as in \eqref{eq:monotonie}, we have for
every $a \in A$,
\begin{eqnarray*}
&& \bigl| u_a^n\bigl(t_k, \xb^1
\bigr) - u_a^n\bigl(t_k, \xb^2
\bigr) \bigr|
\\
&&\qquad\le ( 1 - L_{f,y} \Delta t )^{-1} \biggl( 1 +
L_{f,z} \frac{\E[|\Delta M|] \Delta t }{|a|} + L_{f,\xb} \Delta t \biggr)\bigl |
\xb^1 - \xb^2 \bigr|.
\end{eqnarray*}
It follows that for some constant $C$ independent of $n$,
\[
\bigl| u^n\bigl(t_k, \xb^1\bigr) -
u^n\bigl(t_k, \xb^2\bigr)\bigr | \le
L_{k+1} (1 + C \Delta t),
\]
which implies that $u^n$ is Lipschitz in $\xb$ uniformly for all
$(t_k)_{0 \le k \le n}$ and all $n \ge1$.

(ii) By the Lipschitz property of $u^n$, we have immediately
\[
Du_a^n (t_k, \xb) := \E \biggl[
\frac{u(t_{k+1}, M^{t_k, \xb, a})a^{-1} \Delta
M_{t_{k+1}}^{t_k, \xb, a}}{ \Delta t} \biggr],
\]
is uniformly bounded, which implies that $f(t_k, \xb, u_a^n(t_k, \xb),
D u^n_a(t_k, \xb), a)$ is uniformly bounded.
It follows by the expression \eqref{eq:uan} that
\[
\bigl| u^n(t_k, \xb) - u^n\bigl(t_{k+1},
\xb^{t_k}\bigr) \bigr| \le C \sqrt{\Delta t}.
\]
\upqed\end{pf}

\begin{Proposition} \label{prop:Y_tightness}
Let Assumption~\ref{assum:monotonicity} hold. We have the following properties:
\begin{longlist}[(ii)]
\item[(i)] For every $n \ge1$, there is $e^*_n \in E_n$ such that the
solution $(\Yc^{e^*_n}, \Zc^{e^*_n}, \Nc^{e^*_n})$ of~\eqref
{eq:Discret_BSDE}
satisfies $\Yc^{e^*_n}_{t_k} = u^n(t_k, M^{e^*_n})$, $\P^n$-a.s.

\item[(ii)] The sequence $(\Yc^{e^*_n})_{n \ge1}$ is tight, and $(\Zc
^{e^*_n})_{n \ge1}$ is uniformly bounded.
\end{longlist}
\end{Proposition}

\begin{pf} (i) Let $n \ge1$ be fixed, using the continuity of $H^n_k$ in
$a$ and the dominated convergence theorem,
$a \longmapsto u_a^n(t_k, \xb)$ is continuous, where $u_a^n$ is
defined by \eqref{eq:uan}.
Since $A$ is compact, there is always an optimal $a$ for the
maximization problem \eqref{eq:uan}.
It is then enough to use a classical measurable selection theorem to
construct the required optimal $e^*_n \in E_n$.

(ii) Notice that since we assumed that Assumption~\ref
{assum:monotonicity} holds, we can apply Proposition~\ref{prop:truc}.
Therefore, by \eqref{eq:uan} and using \eqref{eq:un_Lip}, it follows
immediately that
\[
\Zc^{e^*_n}_{t_k} = \E^n_{t_k} \biggl[
\frac{ (u^n(t_{k+1},
M^{e^*_n}_{t_{k+1}}) - u^n(t_k, M^{e^*_n}_{t_k})
)(a^{e_n^*}_{t_k})^{-1}\Delta M^{e^*_n}_{t_{k+1}}}{\Delta t_{k+1}} \biggr],
\]
is uniformly bounded.
Further, using the expression \eqref{eq:uan} with direct computation,
we can easily verify that
\[
\E^n_{t_k} \bigl[ \bigl( \Delta\Yc^{e^*_n}_{t_k}
\bigr)^2 \bigr] \le C \Delta t,
\]
for some constant $C$ independent of $n$, which implies, since $\Yc
^{e^*_n}$ is a pure jump process that
\[
\bigl\langle\Yc^{e^*_n}\bigr\rangle_t\leq Ct_k,\qquad
t_{k-1}\leq t\leq t_{k}.
\]
Finally, we notice that the deterministic nondecreasing process
\[
G^n(s):= C \sum_{k=0}^{n-1}
t_{k+1} \1_{[t_k, t_{k+1})}(s),
\]
converges weakly to the deterministic process $s \longmapsto C s$ as $n
\longrightarrow\infty$.
Then it is enough to apply Theorem~2.3 of Jacod, M\'emin and M\'
etivier \cite{JacodMeminMetivier} for the tightness of $(\Yc
^{e^*_n})_{n \ge1}$, where their condition $C1$ holds for the
nondecreasing process $G^n$.
\end{pf}

\begin{Remark}
In the context of BSDE, Ma, Protter, San Martin and Torres~\cite{MaProtter}
gave a similar tightness result for their numerical solutions,
which is also a key step to prove the convergence of their numerical scheme.
\end{Remark}

%
%
%

Finally, we are ready to provide the proof of Theorem~\ref
{theo:weak_convergence} in two steps.

\begin{pf*}{Proof of Theorem~\ref{theo:weak_convergence}} \textit{Part} (i).
Let us consider the BSDE \reff{eq:BSDE} under some probability measure
$\P\in\mathcal P_S$. In this case, we know that $\overline{\F}^\P
=\overline{\F^{W^\P}}^\P$, for some $\P$-Brownian motion $W^\P$ and
thus that thanks to the predictable representation property, we can
write for some $\overline{\F^{W}}^\P$-predictable process $\tilde a$
\[
B_t = \int_0^t\tilde
a_s^{1/2}\,dW^\P_s,\qquad \P\mbox{-a.s.}
\]

We may now always approach the process $\tilde a$ by a sequence
$(\tilde a^p)_{p\geq0}$ of piecewise-constant processes, over a grid
$(t^p_k)_{0\leq k\leq p}$, whose mesh goes to $0$, in the sense that
\begin{eqnarray*}
\tilde a^p_{t_k}&\in&\overline{\Fc^{W^\P}_{t_k}}^\P\qquad
\mbox{for each $0\leq k\leq p$}\quad\mbox{and}
\\
\E^\P \biggl[\int
_0^T\bigl|\tilde a_s^{1/2}-
\bigl(\tilde a_s^p\bigr)^{1/2}\bigr|^2\,ds
\biggr]&\mathop{\longrightarrow}\limits_{p\rightarrow+\infty}&0.
\end{eqnarray*}

Next, since there is {a priori} no reason that the
applications $\omega\longmapsto\tilde a^p_{t_k}(\omega)$ has any
regularity, we further approximate (by classical density arguments) the
random variables $(\tilde a^p_{t_k})_{0\leq k\leq p}$ by
Lipschitz-continuous functionals\break $(\tilde a^{p,n}_{t_k})_{0\leq k\leq
p}$ such that the following convergence holds true:
\[
\mathbb E^\P \bigl[\bigl|\tilde a^p_{t_k}-\tilde
a^{p,n}_{t_k}\bigr|^2 \bigr]\mathop{\longrightarrow}_{n\rightarrow+\infty}0.
\]

Let us finally denote by $a^n_{\cdot}:=\tilde a^{n,n}_{\cdot}.$ For
every $n
\ge1$, let now $(\Om^n, \Fc^n, \P^n)$ be a probability space
containing $n$ independent random variables $(U_k)_{1 \le k \le n}$,
and consider the following discrete-time martingale defined exactly as
in Section~\ref{sec:main_result}, with functions $H^n_k$ satisfying
\reff{eq:Cond_H}:
\[
M^n_0:=0\quad\mbox{and}\quad M^n_{t_{k+1}}:=M^n_{t_k}+H^n_k
\bigl(a^n_{t_k}\bigl(W^{\P,n}_{\cdot}\bigr),U\bigr),
\]
where $W^{\P,n}$ is a discretized version of $W^\P$ defined by
\[
W^{\P,n}_0:=0\quad\mbox{and}\quad W^{\P,n}_{t_{k+1}}:=W_{t_k}^{\P
,n}+
\bigl(a^n_{t_k}\bigr)^{-1/2}\bigl(W^{\P,n}_{\cdot}
\bigr) H^n_k\bigl(a^n_{t_k}
\bigl(W^{\P,n}_{\cdot}\bigr),U\bigr).
\]

Consider now the following BSDE under $\P^n$
\[
y_t^n=\xi\bigl(\widehat{M}^n_{\cdot}\bigr)-\int
_t^T g_n\bigl(s,
\widehat{M}^n_{\cdot},y_{s^-}^n, a^n_s
z_s^n,a^n_s\bigr)\dvtx d\bigl
\langle M^n_s\bigr\rangle-\int_t^Ta^n_sz^n_s.\,dW^{n}_s-N_T^n+N_t^n,
\]
which is clearly in the same form of BSDE \eqref{eq:Discret_BSDE}, and
hence $y_0^n \le Y_0^n$.

We know that $W^{\P,n}$ converges weakly to $W^\P$. Using
Skrorohod theorem and changing the probability space under which we are
working, it is clear with Lemma~\ref{lemm:M_e} that we may assume
without loss of generality that $W^{\P,n}$ actually converges to $W^\P$
strongly in $\mathcal S^2$ (see also Corollary~14 in Briand, Delyon and
M\'emin~\cite{BriandDelyonMemin} for similar arguments). Moreover,
since the filtrations are Brownian filtrations, we know from \cite
{BriandDelyonMemin} (see their Proposition~3) that the corresponding
filtrations also converge.\footnote{In the sense that, if $(\Fc
_t^n)_{0\leq t\leq T}$ denotes the natural filtration of $W^n$ and
$(\Fc
_t)_{0\leq t\leq T}$ that of $W^\P$, then for every $A\in\mathcal F_T$,
$\E[{\mathbf1}_A|\Fc_t^n]$ converges u.c.p. to $\E[{\mathbf
1}_A|\Fc_t]$.}
Then, using the uniform continuity of $g$ in $t$, we can apply Theorem
$12$ in \cite{BriandDelyonMemin} to obtain that
\[
\lim_{n \to\infty} y^{n}_0 =
\Yc_0^{\P}.
\]
Therefore, we get
\[
\liminf_{n \to\infty} Y_0^n \ge \liminf
_{n \to\infty} y_0^{n} = \Yc^{\P}_0,
\]
which implies the first assertion of Theorem~\ref
{theo:weak_convergence}. 
\end{pf*}

To prove the second part of Theorem~\ref{theo:weak_convergence},\vspace*{1pt}
we shall consider the weak limit of the triplet $(M^{e_n^*}, \Yc
^{e_n^*}, a^{e_n^*})_{n \ge1}$ introduced in Proposition~\ref
{prop:Y_tightness}. Let us first introduce the associated canonical space.
For the process $(M^{e_n^*}, \Yc^{e_n^*})$, it is natural to consider
the spaces of all c\`adl\`ag paths on $[0,T]$ equipped with Skorokhod topology
$ C([0,T], \R^d)$ and $ D([0,T], \R)$ (let us refer to Billinsley
\cite{Billingsley} for a presentation of this canonical space).
For $(a^{e^*_n})_{n \ge1}$, we follow Kushner and Dupuis \cite
{KushnerDupuis} in their numerical analysis to use the canonical space
of measure valued processes (see also El Karoui, Huu Nguyen and
Jeanblanc-Picqu{\'e} \cite{ElKaroui_1987}, or El Karoui and Tan \cite
{ElKaroui_Tan} for a more detailed presentation).
More precisely, since $a^{e_n^*}$ take values in compact set $A$, we
define $\M$ as the space of all measures $m$ on $[0,T] \x A$ such that
the marginal distribution of $m$ on $[0,T]$ is the Lebesgue measure.
By disintegration, $m$ can be write as $m(dt, da) = m_t(da) \,dt$, where
every $m_t$ is a probability measure on $A$, which can be viewed as
measure-valued processes.
We then take $\overline\Om:= C([0,T], \R^d) \x D([0,T], \R) \x\M$
as canonical space, with canonical process $(\overline M, \overline\Yc
, \bar m)$ and the canonical filtration $\overline\F= (\overline\Fc
_t)_{0 \le t \le T}$ generated by the canonical process.
For every $\varphi\in C^2_b(\R^d \x\R)$, we define a process on
$\overline\Om$
\[
\Cc^{\varphi}_t(\overline M, \bar m):= \varphi(\overline
M_t) - \int_0^t \int
_A \frac{1}{2} a \dvtx D^2 \varphi(
\overline M_s ) \bar m_s(da, ds),
\]
and another process
\[
\Dc_t (\overline M, \overline\Yc, \bar m) := Y_t +
\int_0^t \int_A f(s,
\overline M_{\cdot}, \overline\Yc_s, a) \bar
m_s(da, ds),
\]
as well as
\[
\Dc_t^n (\overline M, \overline\Yc, \bar m) :=
Y_t + \int_0^t \int
_A f_n(s, \overline M_{\cdot},
\overline\Yc _s, a) \bar m_s(da, ds),
\]
for every $n \ge1$.
Notice that for every fixed $t > 0$, the two random variables $\Cc
^{\varphi}_t$ and $\Dc_t$ are both bounded continuous in $(\overline M,
\overline\Yc, \bar m)$.

\begin{pf*}{Proof of Theorem~\ref{theo:weak_convergence}} \textit{Part} (ii).
Let us take the sequence $(e_n^*)_{n \ge1}$ introduced in Proposition~\ref{prop:Y_tightness},
we denote $(M^{e_n^*}, \Yc^{e_n^*}, a^{e_n^*})$ by $(M^n, \Yc^n, a^n)$
to simplify the presentation.
Then
\[
\limsup_{n \to\infty} \Yc^n_0= \limsup
_{n \to\infty} Y^n_0.
\]

Denote
\[
m^n(dt , da) := \sum_{k=0}^{n-1}
\delta_{a^n_{t_k}} (da) \,dt 1_{t
\in[t_k, t_{k+1})}.
\]
Let $\Pb^n$ denote the law on $\Omb$ induced by $(\widehat{M}^n, \Yc^n,
m^n)$ in probability space $(\Om^n, \Fc^n, \P^n)$,
where $\widehat M^n$ is the linear interpolation of $(M^n_{t_k})_{0
\le k \le n}$.
Since $(\widehat{M}^n, \Yc^n)_{n \ge1}$ is tight (by Proposition~\ref
{prop:Y_tightness} which uses Assumption~\ref{assum:monotonicity}) and
$A$ is a compact set, then $(\Pb^n)_{n \ge1}$ is relatively compact.
Let $\Pb^{\infty}$ be a limit probability measure,
we claim that
%
\begin{equation}
\label{eq:claim1} \Cc^{\varphi}(\overline M, \bar m) \mbox{ and } \Dc(\overline
M, \overline\Yc, \bar m) \mbox{ are both } \Fbb\mbox{-martingales under }
\Pb^{\infty}.
\end{equation}

Let $0\leq s<t\leq T$ and $\Psi: C([0,T], \R^d) \x D([0,T], \R)
\x\M\longrightarrow\R$ be a bounded continuous function which is
$\Fcb_s$-measurable.
Then by the definition of $(M^n, \Yc^n)$ in \eqref{eq:Discret_BSDE},
it is clear that
\[
\mathbb E^{\overline\P^n} \bigl[\Psi (\overline M, \overline \Yc, \bar m ) \bigl(
\Cc_t^{\varphi}(\overline M, \bar m) - \Cc _s^{\varphi
}(
\overline M, \bar m) \bigr) \bigr] = 0
\]
and
\[
\mathbb E^{\overline\P^n} \bigl[\Psi ( \overline M,\overline \Yc, \bar m )
 \bigl(
\Dc_t^n(\overline M, \overline\Yc, \bar m) - \Dc
_s^n(\overline M, \overline\Yc, \bar m) \bigr) \bigr] =
0.
\]
Since the functionals $\Psi$, $\Cc^{\varphi}_t$ and $\Dc_t$ are all
bounded continuous,
by taking the limit $n \longrightarrow\infty$, it follows that
\[
\mathbb E^{\overline\P^{\infty}} \bigl[\Psi (\overline M, \overline\Yc, \bar m ) \bigl(
\Cc_t^{\varphi}(\overline M, \bar m) - \Cc_s^{\varphi}(
\overline M, \bar m) \bigr) \bigr] = 0
\]
and
\[
\mathbb E^{\overline\P^{\infty}} \bigl[\Psi ( \overline M,\overline\Yc, \bar m ) \bigl(
\Dc_t(\overline M, \overline\Yc, \bar m) - \Dc_s(
\overline M, \overline\Yc, \bar m) \bigr) \bigr] = 0,
\]
which implies claim \eqref{eq:claim1} by the arbitrariness of $\Psi$
and $s \le t$.

It follows that there exists some probability space $(\Om^*, \Fc
^*, \P^*)$ containing the processes $(M^*, \Yc^*, m^{*})$ whose
distribution is $\Pb^{\infty}$.
Let $\F^* = (\Fc^*_t)_{0 \le t \le T}$ be the right-limit of the
filtration generated by $(M^*, \Yc^*, m^*)$, completed under $\P^*$ and
let $a_s^*:=\int_A a m_s^*(da)$ (notice that $a^*$ also takes values
in $A$, since this set is assumed to be convex).
Then $M^*$ is a martingale w.r.t. $\F^*$ with quadratic variation
$\int_0^t a^*_s \,ds$ and $\Dc(M^*, \Yc^*, m^*)$ is a martingale w.r.t. $\F^*$
by claim \eqref{eq:claim1}.
Further, by the convexity of $f$ in $a$, we have
\[
\int_A f(s,\mathbf{x},y,a) m^*_s(da) \geq f
\bigl(s,\mathbf{x},y,a^*_s\bigr).
\]
It follows that $\Yc^{*}_t-\int_0^t f(s,M^*_{\cdot},\Yc_s^{*},a^*_s)\,ds$ is a
bounded $\F^*$-submartingale.

Next, since this is a bounded submartingale, applying Doob--Meyer
decomposition and the orthogonal decomposition for the $\F
^*$-martingales gives us the existence of a $\F^*$-predictable process
$\Zc^*$, a c\`adl\`ag $\F^*$-martingale $\mathcal N^*$, orthogonal to
$M^*$ and a nondecreasing process $\mathcal K^*$ such that
\begin{eqnarray}
\Yc^*_t=\xi-\int_t^Tf
\bigl(s,M^*_{\cdot},\Yc^*_s, a_s^*\bigr)\,ds-\int
_t^T \Zc^*_s \,dM^*_s-
\int_t^T\,d\mathcal N^*_s-\int
_t^T\,d\mathcal K^*_s,\nonumber\\
\eqntext{\P^*
\mbox{-a.s.}}
\end{eqnarray}

Consider now $(\widetilde{\Yc}^*,\widetilde{\Zc}^*,\widetilde
{\mathcal N}^*)$ the unique solution of the following BSDE under~$\P^*$:
%
\begin{eqnarray}
\label{doob_bsde} \widetilde\Yc^*_t=\xi-\int_t^Tf
\bigl(s,M^*_{\cdot},\widetilde\Yc ^*_s,a^*_s\bigr)\,ds-\int
_t^T\widetilde\Zc^*_s\,dM^*_s-
\int_t^T\,d\widetilde{\mathcal
N}^*_s,
\nonumber
\\[-8pt]
\\[-8pt]
\eqntext{\P^*\mbox{-a.s.}}
\end{eqnarray}

We now claim that we necessarily have
%
\begin{equation}
\label{eq:claim2} \mathbb E^{\P^*}\bigl[\widetilde\Yc^*_0\bigr]
\geq\mathbb E^{\P^*}\bigl[\Yc^*_0\bigr].
\end{equation}

This implies that
\[
\limsup_{n \to\infty} Y^n_0=\limsup
_{n \to\infty} \Yc^{e_n}_0 = \mathbb
E^{\P^*}\bigl[\Yc^*_0\bigr]\leq\mathbb E^{\P^*}\bigl[
\widetilde\Yc ^*_0\bigr]\leq \sup_{\P\in\Pc_W}\mathbb
E^\mathbb P\bigl[\Yc_0^\P\bigr] = \sup
_{\P
\in\Pc_S}\Yc_0^\P,
\]
which proves the desired property.

It remains now to prove the claim \eqref{eq:claim2}. It follows
from a classical linearization argument, which we give for
completeness. Using the fact that $f$ is uniformly Lipschitz in $y$, we
may define bounded $\F^*$-progressively measurable process $\lambda$
such that, $\P^*$-a.s.
\[
\delta\Yc^*_t=-\int_t^T
\lambda_s\delta\Yc^*_s \,ds - \int_t^T
\delta\Zc ^*_s\,dM^*_s-\int_t^T\,d
\bigl(\delta\mathcal N^*_s \bigr)+\int_t^T\,d
\mathcal K^*_s,
\]
where
\[
\delta\Yc^*_t:=\widetilde\Yc^*_t-\Yc^*_t,\qquad
\delta\Zc ^*_t:=\widetilde \Zc^*_t-\Zc^*_t,\qquad
\delta\mathcal N_t^*:=\widetilde{\mathcal N}^*_t-\mathcal
N^*_t.
\]

Then denote $\Lambda_t := \exp ( -\int_0^t \lambda_s \,ds
) $.
Applying It\^o's formula to $\Lambda_t\delta\Yc^*_t$ and remembering
that $M^*$ is orthogonal to $\mathcal N^*$ and $\widetilde{\mathcal
N}^*$, we deduce that
\[
\mathbb E^{\P^*}[\delta\Yc_0]=\mathbb E^{\P^*}
\biggl[\int_0^T\Lambda _s\,d\mathcal
K_s \biggr] \geq0,
\]
which completes the proof.
\end{pf*}

\begin{appendix}\label{app}
\section*{Appendix}

We provide here some classical results on BSDEs which are used in the paper.
Let us start by stating a general wellposedness result for BSDEs in an
abstract setting, which will encompass all the cases considered in this paper.

\begin{proposition} \label{prop:sol_GBSDE}
Let $(\Omega_0,\mathcal F,\P)$ be a complete probability space carrying
a square integrable continuous martingale $M$, adapted to a complete
and right-continuous filtration $\F^0:=(\mathcal F_t^0)_{0\leq t\leq
T}$ and a sequence of square-integrable c\`adl\`ag martingales $M^n$
adapted to some filtration $\F^n:=(\mathcal F_t^n)_{0\leq t\leq T}$
which are complete and right-continuous for each $n$. Let $f_0$ and
$f_n$ be functions from $[0,T]\times\Omega_0\times\R\times\R^d$
to $\R
$ and assume furthermore that:
\begin{longlist}[(iii)]
\item[(i)] $\langle M\rangle $ is absolutely continuous with
respect to
the Lebesgue measure, with a density $(a_s)_{0\leq s\leq T}$ taking
values in $A$.
\item[(ii)] There exists a deterministic sequence
$(a_n)_{n\geq
0}$ converging to $0$ such that
\[
\bigl\langle M^n\bigr\rangle _t-\bigl\langle
M^n\bigr\rangle _s\leq C(t-s +a_n)I_d,\qquad
0\leq s\leq t\leq T, \P\mbox{-a.s.},
\]
for some $C>0$.
\item[(iii)] For each $(y,z)$, $f_0(\cdot,M_{\cdot},y,z)$ [resp.,
$f_n(\cdot,M^n_{\cdot},y,z)$] is progressively measurable with respect to
$\F
^0$ (resp., $\F^n$).
\item[(iv)] There is a constant $\mu>0$ such that for each\break 
$n\geq0$ and each $(t,y,y',z,z')$
\begin{eqnarray*}
\bigl|f_0(t,M_{\cdot},y,z)-f_0\bigl(t,M_{\cdot},y',z'
\bigr)\bigr| &\leq&\mu \bigl(\bigl|y-y'\bigr|+\bigl|z-z'\bigr| \bigr),
\\
\bigl|f_n\bigl(t,M^n_{\cdot},y,z\bigr)-f_n
\bigl(t,M^n_{\cdot},y',z'\bigr)\bigr| &\leq&\mu
\bigl(\bigl|y-y'\bigr|+\bigl|z-z'\bigr| \bigr).
\end{eqnarray*}
\item[(v)] For all $(y,z)$, $f_0$ and $f_n$ are continuous in $t$.
\end{longlist}
Then, for $n$ large enough, the following BSDEs under $\P$
%
\setcounter{equation}{0}
\begin{eqnarray}
\label{eq:BSDE_P_W} \Yc_t &=& \xi- \int_t^T
f_0(s, M_{\cdot}, \Yc_s, \Zc_s) \,d
\langle M \rangle _t - \int_t^T
\Zc_s \cdot \,dM_s - \mathcal N_T + \mathcal
N_t,
\\
\qquad\label{eq:BSDE_P_Wn} \Yc_t^n& =& \xi- \int_t^T
f_n\bigl(s, M^n_{\cdot}, \Yc^n_{s^-},
\Zc^n_s\bigr) \,d\bigl\langle M^n \bigr\rangle
_t - \int_t^T \Zc^n_s
\cdot \,dM_s - \mathcal N^n_T + \mathcal
N^n_t,
\end{eqnarray}
where $\mathcal N$ (resp., $\mathcal N^n$) is a c\`adl\`ag
$\F
^0$-martingale (resp., $\F^n$-martingale) orthogonal to $M$
(resp., $M^n$), have a unique solution such that
\begin{eqnarray*}
\mathbb E^\mathbb P \biggl[\sup_{0\leq t\leq T}|\Yc
_t|^2+\int_0^T\bigl|a_s^{{1}/2}
\Zc_s\bigr|^2\,ds+\langle \mathcal N\rangle _T
\biggr]&\leq& C,
\\
\mathbb E^\mathbb P \biggl[\sup_{0\leq t\leq T}\bigl|\Yc
_t^n\bigr|^2+\int_0^T
\Zc_s^n\bigl(\Zc_s^n
\bigr)^T\dvtx d\bigl\langle M^n\bigr\rangle
_s+\bigl\langle \mathcal N^n\bigr\rangle _T
\biggr]&\leq& C
\end{eqnarray*}
for some constant $C>0$ independent of $n$.
\end{proposition}

\begin{pf}
This is actually a direct consequence of the proof of existence via
fixed point arguments in \cite{BriandDelyonMemin}. Indeed, the
assumptions above imply directly that their assumptions $\textup{H1}$, $\textup{H2}$ and
$\textup{H3}$ hold, with the exception that we do not assume that $M^n$
converges to $M$ and that our martingale $M$ can be written as
\[
M_t=\int_0^ta_s^{1/2}\,dW_s,
\]
where $W$ is $(\P,\F^0)$-Brownian motion.

However, by looking carefully at their proofs of Theorem $9$ and
Corollary $10$, it is easy to see that they can be carried out with the
exact same arguments in our setting to obtain the desired results for
the BSDE \reff{eq:BSDE_P_Wn} for $n$ large enough. Moreover, since the
martingale $M$ satisfies their assumption (H1)(ii) with a constant $C
:= \sup_{a \in A} |a|$ and a deterministic sequence $a_n= C |\Delta
_n|$, we can once again follow their proof of existence to obtain
easily that existence, uniqueness and the desired estimates also hold
for \reff{eq:BSDE_P_W}.
\end{pf}


We will now provide a particular robustness result for BSDEs. We go
back to the canonical space $(\Omega,\mathcal F_T)$ and fix a measure
$\P\in\mathcal P_W$. We let $W$ be a $\overline{\F^+}^\P$-Brownian
motion under $\P$, $(a_s)_{0\leq s\leq T}$ be a $\F$-progressively
measurable process and $(a^n)_{n\geq0}$ a sequence of $\F
$-progressively measurable processes such that
%
\begin{equation}
\label{eq:density_approx} \mathbb E^\mathbb P
\biggl[\int_0^T\bigl|{a^n_s}-{a_s}\bigr|^2\,ds
\biggr]\mathop{\longrightarrow}_{n\rightarrow+\infty} 0.
\end{equation}

We next define the following $\overline{\F^+}^\P$-martingales under
$\P$:
\begin{eqnarray*}
M_t&:=&\int_0^ta_s^{1/2}\,dW_s
\quad\mbox{and}\\
\widehat{M}^n_t&:=&\int_0^t
\bigl(a_s^n\bigr)^{1/2}\,dW_s,\qquad
\mathbb P\mbox{-a.s.}
\end{eqnarray*}
Notice that we than have immediately that $\widehat M^n$ converges to
$M$ in the sense that
%
\begin{equation}
\label{eq:S2_approx} \mathbb E^\mathbb P \Bigl[\sup_{0\leq t\leq T}\bigl|M_t-
\widehat {M}^n_t\bigr|^2 \Bigr]
\mathop{\longrightarrow}_{n\rightarrow+\infty
}0.
\end{equation}
We would like to approximate the BSDE
%
\begin{eqnarray}\qquad
\label{eq:BSDE_P_W2} \Yc_t = \xi- \int_t^T
f(s, M_{\cdot}, \Yc_s, \Zc_s, a_s)
\,ds - \int_t^T (a_s)^{1/2}
\Zc_s \cdot \,dW_s - \mathcal N_T + \mathcal
N_t,
\nonumber
\\[-8pt]
\\[-8pt]
\eqntext{\mathbb P\mbox{-a.s.}}
\end{eqnarray}
by the following one for $n\geq0$:
%
\begin{eqnarray}
\label{eq:BSDE_P_W2n} \widehat Y_t^n&=&\xi_n-\int
_t^Tf\bigl(s,\widehat M^n_{\cdot},
\widehat Y_s^n,\widehat Z_s^n,a_s^n
\bigr)\,ds
\nonumber
\\[-8pt]
\\[-8pt]
\nonumber
&&{}-\int_t^T\bigl(a_s^n
\bigr)^{1/2}\widehat Z_s^n\cdot
\,dW_s-\widehat N_T^n+\widehat
N_t^n,\qquad\mathbb P\mbox{-a.s.},
\end{eqnarray}
for some $\mathcal F_T$-measurable random variable $\xi_n$ converging
to $\xi$ in $L^2(\P)$.

\begin{remark}
Notice that existence and uniqueness for these BSDEs are once again
guaranteed by Proposition~\ref{prop:sol_GBSDE}.
\end{remark}

We have the following result, which can be proved using classical
stability arguments for BSDEs. We nonetheless give the proof for completeness.

\begin{proposition} \label{prop:density_approx}
Let Assumptions \ref{assump.gen} hold. Then we have
\[
\mathbb E^{\mathbb P} \biggl[\sup_{0\leq t\leq T}\bigl|\widehat
Y_t^n-\mathcal Y_t\bigr|^2+\int
_0^T\bigl|\widehat Z^n_t-
\mathcal Z_t\bigr|^2\,ds+\bigl\langle \widehat N^{n}-N
\bigr\rangle _T \biggr]\mathop{\longrightarrow}_{n\rightarrow
+\infty
}0.
\]
\end{proposition}

\begin{pf}
Let us apply It\^o's formula to $e^{\eta t}(\widehat Y^n_t-\mathcal
Y_t)^2$, for some constant $\eta$ to be fixed later. We obtain, using
the fact that $\widehat N^n$ and $N$ are orthogonal to $W$
%
\begin{eqnarray}
\label{bsde3}
&&
e^{\eta t}\bigl(\widehat Y^n_t-
\mathcal Y_t\bigr)^2+\int_t^Te^{\eta
s}\bigl|
\bigl(a_s^n\bigr)^{1/2}\widehat
Z^n_s-a_s^{1/2}\mathcal
Z_s\bigr|^2\,ds+\int_t^Te^{\eta s}\,d
\bigl\langle \widehat N^{n}-N\bigr\rangle _s
\nonumber\\
\nonumber
&&\qquad\leq e^{\eta T}|\xi_n-\xi|^2-2\int
_t^Te^{\eta
s}\bigl(\widehat
Y^n_s-\mathcal Y_s\bigr)\\
&&\hspace*{102pt}\qquad{}\times \bigl(f\bigl(s,
\widehat M^n_{\cdot}, \widehat Y_s^n,\widehat
Z_s^n,a_s^n\bigr)-f(s,
M_{\cdot}, \Yc_s, \Zc_s, a_s)
\bigr)\,ds
\\
\nonumber
&&\qquad\quad{}-\eta\int_t^Te^{\eta s}\bigl|
\widehat Y^n_s-\mathcal Y_s\bigr|^2\,ds-2
\int_t^T\bigl(\widehat Y^n_s-
\mathcal Y_s\bigr) \bigl(\bigl(a_s^n
\bigr)^{1/2}\widehat Z^n_s-a_s^{1/2}
\mathcal Z_s \bigr)\cdot \,dW_s
\\
&&\qquad\quad{}-\int_t^Te^{\eta s}\bigl(\widehat
Y_{s^-}^n-\mathcal Y_{s^-}\bigr)\,d\bigl(\widehat
N^n_s-N_s\bigr).\nonumber
\end{eqnarray}

Next, using the uniform continuity of $f$ in $u$ and its Lipschitz
continuity in $(\xb,y,z)$, we have for some modulus of continuity
$\rho
$ and using the trivial inequality $ab\leq\varepsilon a^2+\frac{1}\varepsilon
b^2$ for any $\eps>0$
%
\begin{eqnarray}
\label{bsde4}
\nonumber
&&\biggl|\int_t^Te^{\eta s}
\bigl(\widehat Y^n_s-\mathcal Y_s\bigr)
\bigl(f\bigl(s,\widehat M^n_{\cdot}, \widehat Y_s^n,
\widehat Z_s^n,a_s^n\bigr)-f(s,
M_{\cdot }, \Yc_s, \Zc_s, a_s)
\bigr)\,ds\biggr|
\\
\nonumber
&&\qquad\leq\int_t^Te^{\eta s}\bigl|
\widehat Y^n_s-\mathcal Y_s\bigr|\bigl|f\bigl(s,
\widehat M^n_{\cdot}, \widehat Y_s^n,\widehat
Z_s^n,a_s^n\bigr)-f\bigl(s,
M_{\cdot }, \widehat Y^n_s, \widehat
Z^n_s, a_s\bigr)\bigr|\,ds
\\
&&\qquad\quad{}+\int_t^Te^{\eta s}\bigl|\widehat
Y^n_s-\mathcal Y_s\bigr|\bigl|f\bigl(s,
M_{\cdot}, \widehat Y^n_s, \widehat
Z^n_s, a_s\bigr)-f(s, M_{\cdot},
\Yc _s, \Zc_s, a_s)\bigr|\,ds
\nonumber
\\[-8pt]
\\[-8pt]
\nonumber
&&\qquad\leq C \biggl(\bigl\|\widehat M^n-M\bigr\|_{T}^2+
\int_t^T\rho ^2
\bigl(a^n_s-a_s \bigr)\,ds \biggr)\\
&&\qquad\quad{} + \biggl(C+
\frac{1}\varepsilon \biggr)\int_t^Te^{\eta s}\bigl|
\widehat Y^n_s-\mathcal Y_s\bigr|^2\,ds\nonumber
\\
&&\qquad\quad{}+\varepsilon\int_t^Te^{\eta s}\bigl|\widehat
Z^n_s-\mathcal Z_s\bigr|^2\,ds.\nonumber
\end{eqnarray}

Using the fact that $a^n$ and $a$ are uniformly bounded, if we take
the expectation in \reff{bsde3} and use the estimate \reff{bsde4}, we
obtain by choosing $\eta$ large enough and $\varepsilon<1$
\begin{eqnarray*}
&&\mathbb E^\mathbb P \biggl[\bigl|\widehat Y^n_t-
\mathcal Y_t\bigr|^2+\int_0^T\bigl|
\widehat Z^n_s-\mathcal Z_s\bigr|^2\,ds+
\bigl\langle \widehat N^{n}-N\bigr\rangle _T \biggr]
\\
&&\qquad\leq C\mathbb E^\mathbb P \biggl[|\xi_n-\xi|^2+
\bigl\|\widehat M^n-M\bigr\| _{T}^2+\int
_0^T\rho^2 \bigl(a^n_s-a_s
\bigr)\,ds \biggr].
\end{eqnarray*}

By the dominated convergence theorem and using the fact that $\xi
_n$ converges to $\xi$ and $\widehat M^n$ to $M$, the right-hand side
above goes to $0$. Now the proof can be finished by taking the supremum
in $t$ in \reff{bsde3} and using the BDG inequality. Since this part is
classical, we refrain from writing its proof.
\end{pf}
\end{appendix}

\section*{Acknowledgements}

We are grateful to Nizar Touzi, Jianfeng Zhang and Chao Zhou for
fruitful discussions.
We also would like to thank an anonymous referee and an Associate
Editor, whose advices helped to improve an earlier version of the paper.

%





\printaddresses
\end{document}